\pdfoutput=1
\documentclass[preprint]{elsarticle}

\usepackage{graphicx}
\usepackage{amssymb}

\usepackage{geometry}
\usepackage[utf8]{inputenc} 
\usepackage[T1]{fontenc}    
\usepackage{lineno}
\usepackage{subcaption}
\usepackage{booktabs}       

\usepackage{comment}
\usepackage{color,soul}
\usepackage{epsfig}
\usepackage{subcaption}
\usepackage{amssymb}
\usepackage{amsthm}
\usepackage{amsmath}
\usepackage{multicol}
\usepackage{float}
\usepackage{pslatex}
\usepackage{textcomp}
\usepackage{hyperref}
\usepackage[capitalise, noabbrev]{cleveref}
\usepackage{mathtools}
\usepackage{scalefnt}
\usepackage{cuted}
\usepackage{lipsum}
\usepackage[dvipsnames]{xcolor}
\usepackage[linesnumbered,ruled]{algorithm2e}
\usepackage[printonlyused, withpage]{acronym}
\usepackage{tabularx,array}
\usepackage{makecell}
\usepackage{multirow}
\usepackage{array}
\usepackage{booktabs}
\usepackage{caption}
\usepackage{url}
\usepackage{bm}

\usepackage{todonotes}
\usepackage{mfirstuc}
\usepackage{tikz}
\usepackage{pgfplots}
\usepackage[english]{babel}
\usepackage{xspace}

\pgfplotsset{compat=newest}
\usepgfplotslibrary{groupplots}
\usepgfplotslibrary{dateplot}
\usetikzlibrary{arrows, positioning, shapes.geometric, automata, backgrounds, fit}

\newtheorem{theorem}{Theorem}[section]

\newcommand{\ifdeb}[2]{\ifdefined\debug{#2}\else{#1}\fi}

\newcommand{\mmark}[1]{{\begin{equation} abc \end{equation}}}

\mathchardef\mhyphen="2D

\newcommand{\nonnegintsset}{\mathbb{Z}_{\ge 0}}

\newcommand{\nphard}{{\fontencoding{T1}\fontfamily{cmss}\selectfont NP}-hard\xspace}

\newcommand{\cit}[1]{\ifdeb{\cite{#1}}{Nobel}}
\newcommand{\jobsset}{\ensuremath{J}\xspace}
\newcommand{\jobssetrec}{\ensuremath{\jobsset}\xspace}
\newcommand{\pmax}{\ensuremath{l^{EDD}(\jobsset)}\xspace}
\newcommand{\dmin}{\ensuremath{l^{SPT}(\jobsset)}\xspace}
\newcommand{\authoretal}[1]{#1 \textit{ et al.}}
\newcommand{\authoretalcite}[2]{#1\textit{ et al.} \cit{#2}}
\newcommand{\n}{\ensuremath{n}\xspace}
\newcommand{\rdd}{\ensuremath{\mathit{rdd}}\xspace}
\newcommand{\tf}{\ensuremath{\mathit{tf}}\xspace}
\newcommand{\maxproc}{\ensuremath{p_{max}}\xspace}

\newcommand{\getss}{$\leftarrow$\xspace}

\newcommand{\powten}[1]{\ensuremath{\cdot 10^{#1}}\xspace}

\newcommand{\job}{\ensuremath{j}\xspace}

\newcommand{\numjobs}{\n\xspace}

\newcommand{\verteq}{\rotatebox{90}{$\,=$}}
\newcommand{\equalto}[2]{\underset{\scriptstyle\overset{\mkern4mu\verteq}{#2}}{#1}}

\newcommand{\proctime}[1]{p_{#1}}

\newcommand{\duedate}[1]{d_{#1}}

\newcommand{\tardiness}[1]{T_{#1}}

\newcommand{\position}{\ensuremath{k}\xspace}

\newcommand{\decabstraction}{\ensuremath{\circ}\xspace}

\newcommand{\positionsseth}[1]{\ensuremath{K^{#1}}\xspace}
\newcommand{\positionsset}{\ensuremath{\positionsseth{\decabstraction}(\jobsset)}\xspace}

\newcommand{\positionssetedd}{\ensuremath{\positionsseth{EDD}}(\jobsset)\xspace}
\newcommand{\filpositionssetedd}{\ensuremath{\overline{\positionsseth{EDD}}(\jobsset)}\xspace}
\newcommand{\positionssetspt}{\ensuremath{\positionsseth{SPT}(\jobsset)}\xspace}
\newcommand{\filpositionssetspt}{\ensuremath{\overline{\positionsseth{SPT}}(\jobsset)}\xspace}

\newcommand{\longestjobposition}{\ensuremath{\position^{*}}\xspace}

\ExplSyntaxOn
\DeclareExpandableDocumentCommand{\IfNoValueOrEmptyTF}{mmm}
 {
  \IfNoValueTF{#1}{#2}
   {
    \tl_if_empty:nTF {#1} {#2} {#3}
   }
 }
\ExplSyntaxOff

\DeclareDocumentCommand \pine { o o } {\IfNoValueOrEmptyTF{#1}{}{#2}}
\DeclareDocumentCommand \ifete { o o o } {\IfNoValueOrEmptyTF{#1}{#2}{#3}}

\DeclareDocumentCommand \sequence { o } {
    \IfNoValueTF {#1} {
      \ensuremath{\pi}
    }{
      \ensuremath{\pi_{#1}}
    }
}

\newcommand{\optsequence}[1]{\ensuremath{\pi^*{#1}}\xspace}

\newcommand{\eddsequence}[1]{\ensuremath{\pi^{EDD}_{#1}}\xspace}

\newcommand{\before}[1]{\ensuremath{\sequence[P]}\xspace}

\newcommand{\after}[1]{\ensuremath{\sequence[F]}\xspace}

\newcommand{\precprobh}[2]{\ensuremath{P^{#1}\pine[#2][({#2})]\xspace}}

\newcommand{\follprobh}[2]{\ensuremath{F^{#1}\pine[#2][({#2})]\xspace}}

\newcommand{\precprob}[1]{\precprobh{\decabstraction}{#1}\xspace}

\newcommand{\follprob}[1]{\follprobh{\decabstraction}{#1}\xspace}

\newcommand{\precprobstar}[1]{\ensuremath{\precprob{#1}}}

\newcommand{\jed}[1]{\ensuremath{~\mathrm{#1}}}

\newcommand{\follprobstar}[1]{\ensuremath{\follprob{#1}}}

\newcommand{\precprobedd}[1]{\precprobh{EDD}{#1}}

\newcommand{\follprobedd}[1]{\follprobh{EDD}{#1}}

\newcommand{\precprobspt}[1]{\precprobh{SPT}{#1}}

\newcommand{\follprobspt}[1]{\follprobh{SPT}{#1}}

\newcommand{\acrotable}[1]{\acs{#1} & \acl{#1}\\}

\newcommand{\obj}[1]{\ensuremath{\objposition{#1}{}}\xspace}
\newcommand{\objposition}[2]{\ensuremath{T\pine[#1][_{#2}]\pine[#1][\left(#1\right)]}\xspace}
\newcommand{\objopt}[1]{\ensuremath{T^*\left(#1\right)}}

\newcommand{\predobj}[1]{\ensuremath{\predobjposition{#1}{}}\xspace}
\newcommand{\predobjposition}[2]{\ensuremath{\widehat{\obj{}}\pine[#1][_{#2}]\pine[#1][\left(#1\right)]}\xspace}

\newcommand{\inp}[1]{\ensuremath{\ifete[#1][\bm{X}][\bm{x}_{#1}]\xspace}}
\newcommand{\out}[1]{\ensuremath{y\pine[#1][_{#1}]}\xspace}

\newcommand{\NBR}{\textsc{nbr}\xspace}
\newcommand{\PSK}{\textsc{psk}\xspace}

\ifxetex
  \usepackage{catchfile}
  \newcommand\getenv[2][]{%
    \immediate\write18{kpsewhich --var-value #2 > \jobname.tmp}%
    \newcommand{\temp}{\input{\jobname.tmp}}%
    \if\relax\detokenize{#1}\relax\temp\else\let#1\temp\fi}
\else
  \ifluatex
    \newcommand\getenv[2][]{%
      \edef\temp{\directlua{tex.sprint(
        kpse.var_value("\luatexluaescapestring{#2}") or "" ) }}%
      \if\relax\detokenize{#1}\relax\temp\else\let#1\temp\fi}
  \else
    \usepackage{catchfile}
    \newcommand{\getenv}[2][]{%
      \CatchFileEdef{\temp}{"|kpsewhich --var-value #2"}{\endlinechar=-1}%
      \if\relax\detokenize{#1}\relax\temp\else\let#1\temp\fi}
  \fi
\fi

\newcommand{\drgr}[1]{\IfEq{\LOCAL}{1}{\input{#1}}{\includegraphics[width=0.9\textwidth]{#1}}}

\newcommand{\gapdeffinition}[2]{\ensuremath{
\frac{\obj{\jobsset{#2}, #1}-\objopt{\jobsset{#2}}}
{\obj{\jobsset{#2}, #1}}
}}

\newcommand{\updatemb}{}

\newcommand{\updatembc}[1]{}

\newcommand{\NLP}{\ac{nlp}\xspace}
\newcommand{\LSTM}{\ac{lstm}\xspace}
\newcommand{\GRU}{\ac{gru}\xspace}

\newcommand{\ML}{\ac{ml}\xspace}

\newcommand{\RL}{\ac{rl}\xspace}

\newcommand{\DHS}[1]{\ac{dhs}\pine[{#1}][\kern 0.06em\kern 0.06em\textsc{(\lowercase{{#1}})}]\xspace}

\newcommand{\TTBR}[1]{\ac{ttbr}\pine[{#1}][(\kern 0.09em{#1}~s)]\xspace}
\newcommand{\SMTTP}{$1||\sum T_j$\xspace}
\newcommand{\TSP}{\ac{tsp}\xspace}
\newcommand{\edd}{\ac{edd}\xspace}
\newcommand{\spt}{\ac{spt}\xspace}
\newcommand{\shorter}{\ensuremath{shorter}\xspace}
\newcommand{\OR}{\ac{or}\xspace}
\newcommand{\CSP}{\ac{csp}\xspace}

\newcommand{\gensolve}{\emph{Generate $\&$ Solve}\xspace}
\newcommand{\gensubinstance}{\emph{Subproblem generator}\xspace}

\newcommand{\SOTA}{state-of-the-art\xspace}

\newcommand{\normsumproc}{\textsc{sumproc}\xspace}
\newcommand{\gapeddinv}{\textsc{gapeddinv}\xspace}

\newcommand{\jsplit}{\ensuremath{l^{\decabstraction}(\jobsset)}\xspace}

\acrodef{nlp}[NLP]{\emph{natural language processing}}
\acrodef{lstm}[LSTM]{\emph{Long Short-Term Memory}}
\acrodef{gru}[GRU]{\emph{Gated Recurrent Unit}}
\acrodef{rnn}[RNN]{\emph{recurrent neural network}}
\acrodef{ml}[ML]{machine learning}
\acrodef{ai}[AI]{\emph{artificial intelligence}}
\acrodef{rl}[RL]{\emph{reinforcement learning}}

\acrodef{ga}[GA]{\emph{Genetic Algorithm}}
\acrodef{dhs}[\textsc{horda}]{\emph{Heuristic Optimizer using Regression-based Decomposition Algorithm}}
\acrodef{normSort}[\ensuremath{PSF}]{\emph{Preprocessing Sort Function}}
\acrodef{normConst}[\ensuremath{PC}]{\emph{Preprocessing Constant}}
\acrodef{normCritFunc}[\ensuremath{PCF}]{\emph{Preprocessing Criterion Function}}
\acrodef{normFeaFunc}[\ensuremath{PFF}]{\emph{Preprocessing Feature Function}}
\acrodef{ttbr}[\textsc{ttbm}]{Total Tardiness Branch-and-Merge Algorithm}
\acrodef{smttp}[SMTTP]{Single Machine Total Tardiness Problem}
\acrodef{tsp}[TSP]{Traveling Salesman Problem}
\acrodef{edd}[EDD]{earliest due date}
\acrodef{spt}[SPT]{shortest processing time}
\acrodef{or}[\textsc{OR}]{operations research}
\acrodef{csp}[\textsc{CSP}]{Constraint Satisfaction Problem}

\acrodef{mdd}[\textsc{mdd}]{modified due date rule}

\acrodef{aps}[\textsc{APS}]{Advanced Planning and Scheduling}

\journal{European Journal of Operational Research}

\begin{document}

\newcommand{\LOCAL}{1}

\begin{frontmatter}


\title{Deep learning-driven scheduling algorithm for a single machine problem minimizing the total tardiness}


\acresetall 



\author[1,2]{Michal Bou{\v s}ka}
\ead{bouskmi2@cvut.cz}
\author[1]{P\v remysl \v S\r ucha\corref{cor1}}
\ead{suchap@cvut.cz}
\author[1,2]{Anton{\' i}n Nov{\' a}k}
\ead{antonin.novak@cvut.cz}
\author[1]{Zden{\v e}k Hanz{\' a}lek}
\ead{zdenek.hanzalek@cvut.cz}

\address[1]{Czech Institute of Informatics, Robotics and Cybernetics, Czech Technical University in Prague,\\ Jugosl\'{a}vsk\'{y}ch partyz\'{a}n\r{u} 1580/3, Prague,
Czech republic}
\address[2]{Czech Technical University in Prague, 
Faculty of Electrical Engineering, 
Department of Control Engineering, 
Karlovo~n\'{a}m\v{e}st\'{i} 13, Prague, Czech republic}

\cortext[cor1]{Corresponding author}

\begin{abstract}
In this paper, we investigate the use of the deep learning method for solving a well-known \nphard{} single machine scheduling problem with the objective of minimizing the total tardiness.
    We propose a deep neural network that acts as a polynomial-time estimator of the criterion value used in a single-pass scheduling algorithm based on Lawler\textquotesingle s decomposition and symmetric decomposition proposed by Della Croce et al.
    Essentially, the neural network guides the algorithm by estimating the best splitting of the problem into subproblems.
    The paper also describes a new method for generating the training data set, which speeds up the training dataset generation and reduces the average optimality gap of solutions.
    The experimental results show that our machine learning-driven approach can efficiently generalize information from the training phase to significantly larger instances.
    Even though the instances used in the training phase have from $75$ to $100$ jobs, the average optimality gap on instances with up to 800 jobs is 0.26\%, which is almost five times less than the gap of the state-of-the-art heuristic.
\end{abstract}


\end{frontmatter}


\clearpage
\section{Introduction}
The classical approaches for solving combinatorial problems have several undesirable properties.
First of all, there is a lack of systematic methods that improve the performance of algorithms on unseen instances by gathering the experience from the instances solved in the past.
Therefore, all the information obtained during the past runs of an algorithm is neglected when a new instance is encountered.
A good example is the branch-and-bound-and-remember method~\cite{MORRISON2014403,shang2021}, where the algorithm remembers the information derived during an instance solving, but the information is forgotten as soon as the instance is solved.
Second, the development of efficient heuristic rules requires a substantial amount of time devoted to its design and testing.
This process is tedious and requires a skilled human professional to fine-tune the heuristic's parameters.
A typical example of this feature is genetic algorithms having many parameters for selection, cross-over, mutation, and other operators.

The apparent response to the above challenges is utilizing the existing data.
However, the main obstacle to the successful application of machine learning to enhance algorithms for combinatorial problems remains. It can be formulated as the following fundamental question---\textit{is it possible to extract any useful information from the solved instances and use it efficiently to accelerate solving of an unseen instance?} 

{\updatemb{}This paper addresses a classical \nphard{} single machine total tardiness scheduling problem (\SMTTP), i.e., the problem given by a set of jobs that need to be scheduled on a single machine such that total violation of due-dates is minimized.
Specifically, we investigate the use of deep learning~\cite{krizhevsky2012imagenet}, to guide the solution space exploration of \SMTTP instances.
The presented approach extracts specific information from already solved instances, i.e., parameters of the instance and the optimal value of the objective function.
This information is used as a training data set.
Furthermore, the paper describes a deep neural network that is trained using the training data set. Then the network can predict the optimal value of the objective function for other \SMTTP instances.
Unlike some existing works addressing the use of \ML to solve combinatorial problems (for example, \cit{vinyals2015}), our approach does not rely solely on machine learning, but we combine it with the approaches known from \OR domain.
The described scheduling algorithm shows the way the deep neural network can be combined with classical decomposition schemes \cit{lawler1977, DellaCroce1998} to achieve a fast and efficient solution space exploration.}
The experiments show that our heuristic algorithm outperforms the existing approaches on the standard benchmark data set.
Apart from that, we address the question of how to generate the training data set for the deep neural network. 
A straightforward approach would require solving hundreds of thousands of \nphard{} problems which could take many days.
We show that for problem \SMTTP, there is a much more elegant way that requires only a fraction of that time.

The contributions of this paper can be summarized as follows.
We (i) propose an innovative heuristic algorithm integrating the \ML and \OR approaches;
(ii) improve the process of generating the training data, which leads to faster training and smaller error of our method;
(iii) provide an analysis of deep neural network hyperparameters' impact on the solution's quality, and;
(iv) show that the proposed approach outperforms the state-of-the-art algorithms on the standard benchmark instances.

The rest of the paper is structured as follows. In \autoref{sec:related-work}, we present a review of the literature covering \SMTTP and the use of \ML for solving combinatorial problems.
The studied problem is formally introduced in the subsequent section.
\autoref{sec:dec} describes our approach integrating the \ML into a decomposition-based approach and analyzes its time complexity.
We present results for standard benchmark instances in \autoref{sec:experiments}.
The conclusion is drawn in \autoref{sec:conclusion}, and lists of notations and abbreviations are provided in the \nameref{sec:appendix}.

\section{Related Work}
\label{sec:related-work}
The first part of the related work focuses on the current approaches to solve \SMTTP. 
This part extends the survey by Koulamas \cite{koulamas2010}.
In the second part, we concentrate on existing works exploiting \ML for solving combinatorial problems.

\subsection{Single Machine Total Tardiness Problems}

    In 1977 it was shown by Lawler \cit{lawler1977} that the weighted single machine total tardiness problem is \nphard{}.
    However, it took more than a decade to prove that the unweighted variant of this problem is weakly \nphard{} \cit{du1990}. 
    Lawler \cit{lawler1977} proposed a pseudo-polynomial (in the sum of processing times) algorithm for solving \SMTTP. 
    The algorithm is based on a decomposition of the problem into subproblems. 
    The decomposition firstly sorts the jobs in \edd order.
    Subsequently, it selects the job with the maximum processing time and tries to assign it to the current position and to all the following positions in the \edd sequence.
    For each position, two subproblems are generated.
    The first subproblem contains all the jobs preceding the job with the maximum processing time.
    The second subproblem contains all the jobs following the job with the maximum processing time.
    Besides that, Lawler introduced rules for filtering the candidate positions of the job with the maximum processing time. 
    This algorithm can solve instances with up to one hundred jobs.
    \authoretal{Della Croce} \cit{DellaCroce1998} proposed the \spt decomposition that selects the job with the minimal due date and tries to assign it to every position preceding its original position in the \spt order.
    Similarly to Lawler's decomposition, two subproblems are generated where the first subproblem contains all the jobs preceding the job with the minimal due date, and the second subproblem contains all the jobs following the job with the minimal due date.
    \authoretal{Della Croce} combined both \edd and \spt decompositions together.
    Their algorithm is able to solve instances with up to 150 jobs.
    {
    \updatemb{}
    \authoretal{Szwarc} \cit{szwarc1999} integrated the double decomposition from \cit{DellaCroce1998} and a Split rule \cit{szwarc1996}. Their algorithm solved instances with up to 300 jobs. The same authors further improved the algorithm using paradoxes associated with the problem \cite{szwarc2001}.  
    This algorithm was the state-of-the-art method for a long time, with the ability to solve instances with up to 500 jobs.}
    

    {\updatemb{}Recent papers by \authoretal{Shang}~\cite{Shang2017} and \authoretal{Garraffa} \cit{Garraffa2018} proposed a branch-and-merge algorithm that avoids the solution of equivalent sub-instances in the branching tree. The algorithm uses so-called memorization, i.e., a technique that memorizes the solution of solved sub-problems so that when that sub-problem is reencountered, its solution is retrieved directly from memory instead of solving it again. 
    The authors shown that the algorithm run time converges to $\mathcal{O}^*(2^{\n})$, i.e., the run time is limited by $2^{\n}$ while polynomial factors are omitted.
    The same authors have shown that memorization during the solution space exploration is also efficient for other problems, e.g., $1|r_j|\sum C_j$ and $1|\Tilde{d_j}| \sum w_j C_j$ \cite{shang2021}.
    Nowadays, the algorithm published by Shang, T'Kindt and Della Croce~\cite{shang2021} is the fastest known exact algorithm for \SMTTP, able to solve instances with up to 1200 jobs. In this paper, we denote this algorithm as \TTBR{}.}

    Exact algorithms, such as the ones mentioned above, have very large computation times, while the optimal solution is rarely needed in practice~\cite{WASZ}. 
    Hence, heuristic algorithms are often more practical.
    The existing heuristics algorithms can be categorized into the following three major groups.
    
    The first group of heuristics consists of list scheduling algorithms that create a job order and schedule the jobs according to this order. There are various methods for creating a job order.
    The easiest one is to sort jobs by the Earliest Due Date rule (\edd).
    A more efficient algorithm, called \NBR{}, was proposed in paper \cit{holsenback1992}.
    It is a local search constructive heuristic that starts with job set \jobsset sorted by \edd and constructs the schedule from the end by swapping two jobs by a hand-designed rule.
    \authoretalcite{Panwalkar}{panwalkar1993} proposed an alternative constructive local search heuristic called \PSK{}.
    Russel and Holsenback \cit{russell1997b} compared \PSK{} and \NBR{} heuristics and concluded that neither heuristic is inferior to another one. However, \NBR{} finds a better solution in more cases. 
    
    Heuristics in the second group are based on Lawler's decomposition rule~\cit{lawler1977}. 
    In this case, the heuristic evaluates each node of the search tree, and the most promising node is expanded. This heuristic approach is evaluated in \cit{potts1991} with an \edd heuristic as a guide for the search.
    
    The third group of heuristics contains metaheuristics. 
    Papers \cit{potts1991,antony1996,ben1996} present the simulated annealing algorithm for \SMTTP.
    The same problem is solved in \cit{dimopoulos1999,suer2012} by a genetic algorithm while the authors of \cit{bauer1999,cheng2009} assumed ant colony optimization.
    All the reported results in the previous studies are for instances with up to $100$ jobs.
    However, these instances are solvable by the current state-of-the-art exact algorithm in a fraction of a second.
    
\subsection{Use of Machine Learning in Algorithms for Combinatorial Optimization Problem}
    The integration of \ML into algorithms for solving combinatorial optimization problems has several difficulties.
    {\updatemb{}
    First, the instances of scheduling problems naturally appear in different sizes, e.g.,  with a variable number of jobs. 
    In opposite to this, the majority of \ML models are often designed with a fixed size of the input feature vector and the output vector.
    }
    This issue can be addressed by recurrent networks and, more recently, by encoder-decoder type of architectures.
    Vinyals~\cit{vinyals2015} applied an architecture called Pointer Network that, given a set of graph nodes, outputs a solution as a permutation of these nodes.
    The authors applied the Pointer Network to \TSP with up to 20 nodes; however, this approach for \TSP is still not competitive with the best classical solvers such as Concorde~\cit{applegate2006concorde} that can find optimal solutions to instances with $80,000$ nodes.
    Moreover, the Pointer Network output needs to be corrected by the beam-search procedure, which underlines the weaknesses of this end-to-end approach.
    Pointer Network has achieved an optimality gap of around 1\% for instances with 20 nodes after performing the beam search.
    
    The second difficulty with using \ML models for solving combinatorial problems lies in the acquisition of training data.
    Obtaining a single label for a training instance usually requires solving a problem of the same complexity as the original problem itself, while \ML usually requires millions of training samples.
    This issue can be addressed by the reinforcement learning paradigm.
    \authoretal{Deudon}~\cite{deudon2018learning} used encoder-decoder architecture trained with REINFORCE algorithm to solve 2D Euclidean \TSP with up to 100 nodes.
    It is shown that (i) repetitive sampling from the network is needed, (ii) applying a well-known 2-opt heuristic on the results still improves the solution of the network, and (iii) both the quality and run times are worse than classical exact solvers.
    A similar approach, used to solve TSP, is described in~\cite{kool2018} which, if it is treated as a greedy heuristic, beats simple heuristics such as Christofides algorithm on small instances.
    To be competitive with a relevant baseline algorithm such as Lin-Kernighan heuristics \cite{lin1973}, they perform repeated sampling from the model and output the best solution.
    Moreover, they do not directly compare their approach with state-of-the-art classical algorithms while admitting that general-purpose Integer Programming solver Gurobi solves their largest instances optimally within 1.5\jed{s}.
    
    Reinforcement learning was also used to solve other combinatorial problems.
    For example, \authoretal{Khalil} \cite{khalil2017learning} presented an approach for learning greedy algorithms over graph structures. 
    The authors show that their S2V-DQN model can obtain competitive results on MAX-CUT and Minimum Vertex Cover problems. For \TSP, S2V-DQN performs about the same as 2-opt heuristics. 
    Unfortunately, the authors do not compare running times with Concorde solver.
    Interesting results for graph coloring were introduced by \authoretal{Huang} \cite{huang2019}.
    Huang~\textit{et~al.} proposed a \RL heuristic with a neural network able to outperform the state-of-the-art heuristic by 1-2\%, when trained on the same type of graph as the one used during the evaluation.
    \authoretal{Abe}~\cite{abe2019} presented an \RL approach for Minimum Vertex Cover and MAX-CUT. 
    For Minimum Vertex Cover problem, they have up to 10\% better solutions than the 2-approximation algorithm. 
    For MAX-CUT problem, they are not able to outperform the heuristic of Laguna \cite{laguna2009}.
    More details can be found in the survey by Mazyavkina~\textit{et~al.}~\cite{mazyavkina2020} addressing the use of \RL approaches for solving combinatorial problems.

    Integration of \ML with scheduling problems has received little attention so far.
    Earlier attempts of integrating neural networks with job-shop scheduling were published in~\cit{zhou1991} and \cit{jain1998}.
    However, their computational results are inferior to the traditional algorithms, or it is not possible to assess their quality.
    Alternative use of \ML in the scheduling domain is focused on the evaluation of criterion functions.
    For example, the authors in~{\cit{vaclavik2016}} addressed a nurse rostering problem and proposed a classifier, implemented as a neural network, able to determine whether a particular change in a solution leads to a better solution or not.
    This classifier is then used in a local search algorithm to filter out solutions having a low chance of improving the criterion function.
    Nevertheless, the approach is sensitive to changes in the problem size, i.e., the length of the planning horizon.
    If the size of the problem is changed, a new neural network must be trained.
    Another method, which does not directly predict a solution to the given instance, is proposed in~\cit{novak2015}.
    In this case, an online \ML technique is integrated into an exact algorithm where it acts as a heuristic.
    Specifically, the authors use regression for predicting the upper bound of a pricing problem in a Branch-and-Price algorithm.
    Correct prediction leads to faster computation of the pricing problem, while incorrect prediction does not affect the optimality of the algorithm.
    This method is not sensitive to the change of the problem size; however, it is designed specifically for the Branch-and-Price approach and cannot be generalized to other approaches.
    {
    \updatemb{}
    The authors of paper \cite{lara2019} use the neural network as hyper-heuristic switching between several known heuristics for the job-shop scheduling problem. Nevertheless, it is hard to assess the benefit of the method since a comparison with existing approaches is not provided.
    Recently, \authoretal{Zhang} \cite{zhang2020} solved the same scheduling problem using end-to-end deep reinforcement learning.
    The authors trained Graph Neural Network to generate a priority dispatching rule. The results show that their approach provides better results compared to simple priority rules like Shortest Processing Time, Most Work Remaining, etc.
    Shu Luo~\cite{luo2020} proposed a Deep Q-Network (DQN) trained by reinforcement learning for the online dynamic flexible job shop scheduling problem.
    The neural network is trained to select one out of six dispatching rules. The selected rule is then applied in each iteration of their algorithm.
    A similar work of~\cite{alicastro2021} applies reinforcement learning to an additive manufacturing machine scheduling problem. The authors describe a reinforcement learning iterated local search meta-heuristic that switches different operators of the local search.

    In the field of \CSP solving, Xu \textit{et al.} \cite{xu2018} presented a neural network estimating the satisfiability of the \CSP.
    The authors assume that the neural network can be integrated into an algorithm for solving \CSP.
    However, Xu \textit{et al.} tested their approach only on instances with up to 128 binary variables.
    \authoretal{Cappart} \cite{cappart2020} trained neural network able to estimate values of variables during of \CSP solving.
    The estimation of the value of a variable with a neural network leads to an earlier finding of the part of the state space with an optimal solution and thus to faster convergence.
    The approach is able to solve the instances of Travelling Salesman Problem with Time Windows with up to $100$ nodes.
    Although the proposed method shows an interesting idea, it should be noted that in the literature, there are classical approaches that solve instances with up to $200$ nodes.

    \authoretal{Nair}~\cite{nair2020} introduced an approach to speed up a Mixed Integer Linear Programming solver with a neural network.
    They present two methods to speed up the solution --- Neural Diving and Neural Branching. 
    Neural Diving focuses on the improvement of the incumbent bound. 
    It generates a partial solution of the instance (i.e., predicts a value only for a subset of variables), which is then fixed. 
    The values of the remaining variables are solved by the solver. 
    Neural Branching is used to select a branching variable in the \textit{branch-and-bound} method.
    It aims to approximate a computationally expensive branching strategy with just a fraction of the computation time.
    By the combination of these two methods, \authoretal{Nair} achieves 1.5 times smaller optimality gap on MIPLIB benchmark set.
    Another different way to speedup MILP solvers is introduced by  \authoretal{Tang} \cite{tang2020}.
    \authoretal{Tang} trained the \RL agent, which learns to generate cutting planes and is shown to outperform the human-designed heuristic used in Gurobi MILP solver.
    Their approach achieves $2$ to $3$ times faster convergence on large instances for packing, production planning, and MAX-CUT problems.
    }

    The contemporary operations research literature has started to focus on machine learning approaches more intensively.
    A recent survey by Bengio~\textit{et al.}~\cite{bengio2020} identifies four main problems of the use of \ML in combinatorial optimization, i.e., modeling, feasibility, scaling, and data generation.
    (i) Bengio~\textit{et al.} argue that unlike, e.g., computer vision, there are no neural network models in the literature that would be suitable for combinatorial problems. 
    (ii) Apart from that, neural networks can be used only as a heuristic. Therefore, their current use is limited for exact approaches as well as for problems where it is difficult to find a feasible solution.
    (iii, iv) The last two problems correlate with the first two paragraphs in this section.
    The authors conclude that the existing approaches are at an early stage of development, but they open new opportunities for research addressing combinatorial optimization algorithms.
    
    {\updatemb{}
    This paper is founded on the idea that the frequent limitation of the existing techniques applying ML to combinatorial problems is the use of end-to-end approaches. Their weakness is that they disregard fundamental properties of the combinatorial problems that have been studied in the literature for decades.
    The view studied in this paper is different, and the proposed solution efficiently combines both the ML and properties of the problem.
    }

\section{Problem Statement}
\label{sec:problem}

    In this paper, we study a single machine scheduling problem defined by a set of jobs $\jobsset = \{1,\dots,\numjobs\}$.
    The machine can process at most one job at a time, and all the jobs are available for processing at time zero.
    The execution of the jobs cannot be interrupted.
    Each job $\job \in \jobsset$ is characterized by processing time $\proctime{\job} \in \mathbb{Z}_{>0}$ and due date $\duedate{\job} \in \nonnegintsset$.
    {\updatemb{}
    A solution to this problem is a \emph{schedule} given by a one-to-one correspondence $\sequence: {\{1,\dots,\numjobs\}} \mapsto {\{1,\dots,\numjobs\}}$ 
    mapping a position in the schedule to a job, i.e., $\sequence[\position] \in \jobsset$ is the job at position $\position$ in schedule $\sequence$.}
    For a scheduled job, the problem defines its tardiness as an indicator measuring how much the job violates the due date.
    Tardiness of job $\sequence[\position]$ in schedule $\sequence$ is defined as $\mathcal{\tardiness}_{\sequence[\position]}(\jobsset) = \max\left(0, \sum_{\position^{\prime} \in \jobsset: \position^{\prime} \leq \position} \proctime{\sequence[\position^{\prime}]} - \duedate{\sequence[\position]}\right)$.
    Then, the total tardiness of schedule \sequence is defined as $\obj{\jobsset,\sequence} = \sum^\numjobs_{\position=1} \mathcal{\tardiness}_{\sequence[\position]}(\jobsset)$.
    The goal of the scheduling problem is to find an optimal schedule \optsequence{} which minimizes the total tardiness $\objopt{\jobsset} = \min_{\sequence \in \Pi} \obj{\jobsset,\sequence}$ where $\Pi$ is the set of all jobs' permutations.
    To ease the readability of the paper, the list of notation and symbols used is provided in Appendix.

    This combinatorial problem is proven to be weakly \nphard{}~\cit{du1990}. Graham's notation~\cit{graham1979} denotes it as $1||\sum \tardiness{\job}$ where 1 indicates that it is a single machine scheduling problem and $\sum \tardiness{\job}$ refers to the objective function, i.e., $\min_{\sequence \in \Pi} \obj{\jobsset,\sequence}$.

\section{Proposed Decomposition Heuristic Algorithm}
\label{sec:dec}
In this section, we introduce \DHS{} for \SMTTP.
The algorithm's name comes from the fact that it uses decomposition controlled by a regressor.
The regressor is realized using a deep neural network (neural network for short in the rest of the paper) approximating the relation between features of instance and \objopt{\jobsset}.

The description of the algorithm is structured as follows.
First of all, we summarize the problem decompositions used in the \DHS{} algorithm.
Second, we describe \DHS{} and show how it effectively combines the well-known properties of \SMTTP and an \ML model.
Next, we proceed by discussing the architecture of the regressor and its integration into \SMTTP decompositions, and we describe training data acquisition, including the training of the neural network.
Finally, we analyze the time complexity of \DHS{} algorithm.
    
In the rest of the paper, we use two definitions to describe the ordering of the job set \jobsset:
\begin{enumerate}
    \item \ac{edd}: if $1 \le \job < \job^{\prime} \le \numjobs$ then either (i) $\duedate{\job} < \duedate{\job^{\prime}}$ or (ii) $\duedate{\job} = \duedate{\job^{\prime}} \, \wedge \,\proctime{\job} \le \proctime{\job^{\prime}}$,
    \item \ac{spt}: if $1 \le \job < \job^{\prime} \le \numjobs$ then either (i) $\proctime{\job} < \proctime{\job^{\prime}}$ or (ii) $\proctime{\job} = \proctime{\job^{\prime}} \, \wedge \,\duedate{\job} \le \duedate{\job^{\prime}}$.
\end{enumerate}
\ac{edd}~(Earliest Due Date) is a sequence of jobs, sorted in non-decreasing order of due dates and \ac{spt}~(Shortest Processing Time) is a sequence of jobs sorted in non-decreasing time of processing times.

\subsection{Problem Decompositions}
    \label{sec:smttpdec}
    \label{sec:dhs:decompositions}    
    Before we describe \DHS{}, it is necessary to outline two decomposition approaches for \SMTTP that are used in our algorithm. 
    The first decomposition is Lawler's decomposition~\cit{lawler1977} which utilizes \ac{edd} order of jobs; therefore, the related notation is denoted by superscript \ac{edd}. The other decomposition, proposed by Della Croce \textit{et al.}~\cit{DellaCroce1998}, is analogous but based on \ac{spt} order of jobs; thus, the related notation is denoted as \ac{spt}.

    Both decompositions exploit the fact that any optimal schedule of \SMTTP can be represented by a permutation of jobs \sequence since the machine is never idle in an optimal schedule.
    Each decomposition~$\decabstraction \in \{\ac{edd},\ac{spt}\}$ defines the splitting job $\jsplit \in \jobsset$, i.e., \pmax for Lawler's decomposition and \dmin for the decomposition proposed by Della Croce \textit{et al.}
    For a position \position of job \jsplit in the schedule, the decomposition splits \jobsset into two subsets.
    The first subset $\precprob{\jobsset, \position}$ represents jobs preceding \jsplit in the schedule, and the second subset $\follprob{\jobsset, \position}$ represents jobs following job \jsplit in the schedule under decomposition~\decabstraction.
    The precise definition of the $\precprob{\jobsset, \position}$ and $\follprob{\jobsset, \position}$ is linked with a particular decomposition, as it is explained below.
    
    The first decomposition, further denoted as \edd decomposition, is based on a theorem proposed by Lawler \cit{lawler1977}.
    \begin{theorem}
        \label{theo:law-dec}
        (Lawler, 1977)
        Suppose jobs \jobsset are ordered in \edd order and the splitting job is $\pmax = \arg \max_{i \in \jobsset}{p_i}$.
        Then, there is some integer \position, $\pmax \le \position \le \numjobs$, such that there exists an optimal sequence $\optsequence{}$ in which the splitting job \pmax is preceded by all jobs \job such that $\job \le \position$, and followed by all jobs \job such that $j > \position$.
    \end{theorem}

    Lawler's decomposition splits jobs into two subsets $\precprobedd{\jobssetrec,\position}$ and $\follprobedd{\jobssetrec,\position}$, which for job set \jobsset and position \position contains jobs $\{1,\dots,\position \} \setminus \{ \pmax \}$ and $\{\position + 1, \dots, \numjobs \}$, respectively.
    Thus, for each eligible position $\position \in \{\pmax, \dots, \numjobs \}$, the problem is decomposed into two subproblems defined by \precprobedd{\jobsset, \position} and \follprobedd{\jobsset, \position} such that job \pmax is neither in \precprobedd{} nor in \follprobedd{}.
    When we denote the set of positions $\{\pmax, \dots, \numjobs \}$ as $\positionssetedd$, then  the optimal total tardiness $\objopt{\jobsset}$ of instance \jobsset can be computed as
    \begin{equation}
        \objopt{\jobsset}  =\min_{\position \in \positionssetedd} Q(\jobsset,\position),
        \label{eq:decompObjective}
    \end{equation}
    where
    \begin{equation}
        \label{eq:exact-criterion}
        \begin{split}
            Q(\jobsset,\position) =
            \objopt{\precprobedd{\jobssetrec,\position}}
            +
            \max\left\{0,
            \sum_{\job \in \precprobedd{\jobssetrec,\position}}
            \proctime{\job}
            + \proctime{\position}
            - \duedate{\position}
            \right\} + 
            \objopt{\follprobedd{\jobssetrec,\position}}\,.
        \end{split}
    \end{equation}
    The optimal solution to the instance is found by recursively selecting position \position with the minimal criterion $Q(\jobsset,\position)$.

    The second decomposition, denoted as \spt decomposition, was proposed by \authoretal{Della Croce} \cit{DellaCroce1998} and is described by the following theorem.
    \begin{theorem}
        \label{theo:spt-dec}
        (Della Croce, 1998)
        Suppose jobs \jobsset are in \spt order and job $\dmin = \arg \min_{i \in \jobsset}{\duedate{i}}$.
        Then there exists an integer \position, $1 \le \position \le \dmin$, such that there exists and optimal sequence $\optsequence{}$ in which the jobs preceding \dmin are uniquely determined as follows: take jobs $\{1, \ldots, \dmin - 1\}$ in \spt order and sort these jobs by the \edd order and select the first $\position - 1$ jobs.
        All other jobs follow the \dmin job.
    \end{theorem}
    \autoref{theo:spt-dec} describes a similar decomposition to \edd decomposition but uses different position set $\positionssetspt = \{1,\ldots, \dmin \}$.
    Furthermore, the set of preceding jobs is denoted as \precprobspt{\jobsset, \position}, while the set of following jobs is \follprobspt{\jobsset, \position}.
    Nevertheless, the basic idea of the decomposition is the same as the one formulated in Equation~(\ref{eq:decompObjective}) for the \edd decomposition.

    The efficiency of both decomposition approaches is significantly influenced by the number of the relevant positions for the splitting job \jsplit, $\circ\in\{\edd,\spt\}$, i.e., cardinalities of $\positionssetedd$ and $\positionssetspt$.
    The size of the position set $\positionsset$ for $\circ\in\{\edd,\spt\}$ can be reduced by filtering rules described in \cit{lawler1977,szwarc1996,szwarc1999}.
    These rules can exclude some positions that provably cannot lead to an optimal solution.
    {
    \updatemb{}The efficiency of the decompositions can be improved by the following rules:
    (i) remove the position from $\positionsset{}$ if the completion time of job $\jsplit$ at position $\position$ is larger than the due date of the following job~\cit{szwarc1999},
    (ii) remove the position from $\positionsset{}$ if the completion time of job $\jsplit$ at position $\position$ is smaller than the due date plus the processing time of the previous job~\cit{szwarc1999}.
    The other rules are based on a similar idea. A detailed explanation, including the proof, can be found in \cite{szwarc1999}. 
    In the rest of the paper, we denote the filtered set of positions by $\filpositionssetedd \subseteq \positionssetedd$, $\filpositionssetspt \subseteq \positionssetspt$, i.e., $\overline{K^{\circ}}(\jobsset) \subseteq \positionsset$ for $\circ\in\{\edd, \spt\}$.      
    }
    
\subsection{Scheduling Algorithm}
    \label{sec:dhs:dhs}
    Even though algorithms, e.g., \cite{Garraffa2018}, that use decompositions described in the previous section, are very efficient, their time complexity grows very quickly with the number of jobs.
    Therefore, we propose a heuristic algorithm, denoted as \DHS{}, which approximates the search of the solution space by \textit{a priori} trained regressor.
    
    \DHS{} is a greedy heuristic that combines the efficiency of \edd and \spt decompositions and \ML. 
    \DHS{} recursively applies one of the decompositions while the position \position of the splitting job \jsplit is determined using the regressor.  It aims to select the best position $k^*$ of the splitting job \jsplit in \filpositionssetedd or \filpositionssetspt without solving the subproblems. Therefore, the regressor estimates values of $T^{*}(\precprob{\jobsset, \position})$ and $T^*(\follprob{\jobsset, \position})$ in \cref{eq:exact-criterion} for every relevant position $\position\in \overline{K^{\circ}}(\jobsset)
$.
    With that, the algorithm selects position $\position^*$ that minimizes the estimate of the objective function.

    \begin{algorithm}[!htb]
    \SetKwInput{kwInit}{Input}
        \DontPrintSemicolon
        \kwInit{a set of jobs \jobsset}
        \KwOut{schedule \sequence{}}
        \SetKwFunction{FDHS}{\DHS{}}
        \SetKwProg{Fn}{Function}{:}{}
        \Fn{\FDHS{\jobssetrec}}{ \label{code:dhs:TTBR1}
     {\updatemb{}\tcc{Small instances are solved using an exact approach.}}
     \If{$|\jobssetrec| \leq 5$}
     { \label{code:dhs:small1}
        \sequence \getss \TTBR{}(\jobsset)\;
        \Return{\sequence} \label{code:dhs:TTBR2}
     }\label{code:dhs:small2}
     \tcc{Determines the splitting job and the position set (for both decompositions).}
     \pmax, \filpositionssetedd \getss genEDDPos(\jobssetrec)\;\label{code:dhs:genedd}
     \dmin, \filpositionssetspt \getss genSPTPos(\jobssetrec)\;\label{code:dhs:genspt}
     {\updatemb{}\tcc{Selection of position set with smaller cardinality.}}
     \If{$|\filpositionssetedd| \le |\filpositionssetspt|$\label{code:dhs:selectstart}}
     {
         $\overline{K^{\circ}}(\jobsset)$
 \getss \filpositionssetedd; \jsplit \getss \pmax\;
     }\Else
     {
        $\overline{K^{\circ}}(\jobsset)$
 \getss \filpositionssetspt; \jsplit \getss \dmin\;
     }\label{code:dhs:selectend}
        
     \tcc{Determine the position of the splitting job using the regressor.}
     \longestjobposition \getss  $\arg \min_{\position \in \overline{K^{\circ}}(\jobsset)}\,(\predobj{\precprobstar{\jobssetrec, \position}} + \max(0, \proctime{\position} - \duedate{\position} + \sum_{\job \in \precprobstar{\jobsset{}, \position}}\proctime{\job}) +  \predobj{\follprobstar{\jobssetrec, \position}})$\;  \label{code:dhs:argmin}
     {\updatemb{}\tcc{Recursively call of \DHS{} for both subproblems.}}
     \before{} \getss \FDHS{\precprobstar{\jobssetrec, \longestjobposition}}\; \label{code:dhs:recursionb}
     \after{} \getss \FDHS{\follprobstar{\jobssetrec, \longestjobposition}}\; \label{code:dhs:recursiona}
     \tcc{join sequences into one}
     \sequence{} \getss (\before{}, \jsplit, \after{})\; \label{code:dhs:merge}
    \Return{\sequence{}} \label{code:dhs:return}
    }
    \caption{Heuristic Optimizer using Regression-based Decomposition Algorithm (\DHS{})}
    \label{code:dhs}
\end{algorithm}

    {\updatemb{} The algorithm is outlined in \autoref{code:dhs}.
    It recursively applies one of the decompositions described in the previous section while splitting input jobs set \jobsset to its subsets \precprobstar{\jobssetrec, \longestjobposition} and \follprobstar{\jobssetrec, \longestjobposition}.
    In the first step (lines \ref{code:dhs:small1}~and~\ref{code:dhs:small2}), the algorithm handles job sets with five or fewer jobs.
    It turned out to be more efficient to run an exact algorithm instead of performing the inference from the regressor on a small set of jobs.
    Subsequently, the algorithm determines the splitting job and set of positions \position for both decompositions, i.e., \pmax, \filpositionssetedd and \dmin, \filpositionssetspt for \edd and \spt decomposition, respectively (lines \ref{code:dhs:genedd}~and~\ref{code:dhs:genspt}).
    To reduce the run time, \DHS{} uses either \edd or \spt decomposition depending on the cardinalities of their position sets (lines \ref{code:dhs:selectstart}--\ref{code:dhs:selectend}). 
    The position set with smaller cardinality is selected and the selected position set and the related splitting job are stored to $\overline{K^{\circ}}(\jobsset)$ and \jsplit, respectively.
    After the selection of the positions set, the algorithm greedily selects position \longestjobposition (line \ref{code:dhs:argmin}) having the minimal estimation of the optimal objective function, i. e., the total tardiness. Thus the position is determined as $\longestjobposition = \arg \min_{\position \in \overline{K^{\circ}}(\jobsset)
} \hat{Q}(\jobsset,\position)$, where $\hat{Q}(\jobsset,\position)$ is the estimate of the objective function for position \position computed as}
    \begin{equation}
        \label{eq:sol:dhs:predicted-criterion:jobsset}
        \begin{split}
            \hat{Q}(\jobsset,\position) =
            \predobj{\precprob{\jobssetrec,\position}} + 
            \max\left\{0,
            \sum_{\job \in \precprob{\jobssetrec,\position}}
            \proctime{\job}
            + \proctime{\position}
            - \duedate{\position}
            \right\} + 
            \predobj{\follprob{\jobssetrec,\position}}\,.
        \end{split}
    \end{equation}
    Estimates \predobj{\precprob{\jobsset, \position}} and \predobj{\follprob{\jobsset, \position}} are computed by the regressor described in the following section.
    Subsequently, the algorithm recursively solves job sets   \precprob{\jobssetrec,\longestjobposition} and \follprob{\jobssetrec,\longestjobposition}.
    Resulting partial sequences are stored as vectors \before{} and \after{}, respectively (lines \ref{code:dhs:recursionb}~and~\ref{code:dhs:recursiona}).
    Finally, the algorithm merges $(\before{}, \jsplit, \after{})$ into one sequence $\sequence{}$, which is returned as the resulting schedule (\autoref{code:dhs:return}).
    
    {
    \updatemb{}Please notice that if the estimates \predobj{\precprobstar{\jobssetrec, \position}} and \predobj{\follprobstar{\jobssetrec, \position}} in Equation (\ref{eq:sol:dhs:predicted-criterion:jobsset}) would be perfect (i.e., \predobj{\precprobstar{\jobssetrec, \position}} = \objopt{\precprob{\jobssetrec,\position}} and \predobj{\follprobstar{\jobssetrec, \position}} = \objopt{\follprob{\jobssetrec,\position}} for every \position $\in$ $\overline{K^{\circ}}(\jobsset))$, then \DHS{} would turn into an exact method. This claim directly follows from theorems \ref{theo:law-dec} and \ref{theo:spt-dec}. Since Algorithm~\ref{code:dhs} enumerates all positions \position that may lead to an optimal solution, then accurate estimates \predobj{\precprobstar{\jobssetrec, \position}} and \predobj{\follprobstar{\jobssetrec, \position}} guarantee finding \position leading to $\optsequence{}$ according to theorems \ref{theo:law-dec} and \ref{theo:spt-dec}. The first iteration of Algorithm~\ref{code:dhs} with perfect estimates finds the optimal value of the objective function, while its recursive application finds $\optsequence{}$. 
    Nevertheless, assuming that {\fontencoding{T1}\fontfamily{cmss}\selectfont{ P}$\neq$\fontencoding{T1}\fontfamily{cmss}\selectfont{NP}\fontencoding{T1}\fontfamily{lmr}\selectfont}, it is impossible to guarantee that estimates are always correct.
    On the other hand, the better the precision of the regressor, the more likely that \DHS{} will find the optimal solution.
    
    We also experimented with other improvements, known in the scheduling domain, like the Split rule proposed in \cit{szwarc1996}. The rule extends the filtering rules and allows to determine the optimal block sequence w.r.t. position \position. However, the tests have shown that in our setting, it slows down \DHS{} and does not improve the quality of the solutions. Thus we do not use it in Algorithm~\ref{code:dhs}.
    }

\subsection{Regressor}
    \label{sec:dec:regressor}

    \DHS{} algorithm utilizes the regressor to estimate $\hat{T}(\jobsset^\prime)$ where $\jobsset^\prime$ is either $\precprob{\jobssetrec,\position} \subset \jobsset$ or $\follprob{\jobssetrec,\position} \subset \jobsset$.
    It comes as no surprise that the quality of the estimation significantly affects the ability of \DHS{} to find an optimal or near-optimal solution.
    The key advantage of \DHS{} is that it is not sensitive to the absolute error of the estimation since the algorithm compares multiple solutions obtained by the same regressor for different $\position \in \overline{K^{\circ}}(\jobsset)$.
    Therefore, the proposed regressor uses a neural network since those have been shown to be successful for problems being sensitive to relative error~\cite{vaclavik2016}.

    \begin{figure}[!htb]
        \centering
        \scalebox{0.85}{
        \begin{tikzpicture}[
            node distance = 0.95cm and 0.95cm,
            process/.style = {rectangle, draw, fill=gray!30,
                            minimum width=1.8cm, minimum height=0.85cm, align=center},
            data/.style = {rectangle, minimum height=0.85cm, align=center},
            arrow/.style = {thick,-stealth},
            bendarrow/.style = {thick,-stealth},
                                ]
            \node (n1)  [data] {$\jobsset{}' \approx \left\{(p_j, d_j)\right\}_{j\in\jobsset^\prime}$};
            \node (n2)  [process,right=of n1]   {$Norm$};
            \node (n3)  [process,right=of n2]   {recurrent\\layer};
            \node (n4)  [process,right=of n3]   {dense\\layer};
            \node (n5)  [process,right=of n4]   {$Norm^{-1}$};
            \node (n6)  [data,right=of n5]   {$\hat{T}(\jobsset') \in \mathbb{R}_{\geq 0}$};
            \draw [arrow] (n1) -- (n2);
            \draw [arrow] (n2) -- (n3)node[midway,below] {$\equalto{\inp{}}{Norm(\left\{(p_j, d_j)\right\}_{j\in\jobsset^\prime})}$};
            \draw [arrow] (n3) -- (n4)node[midway,right] {};
            \draw [arrow] (n4) -- (n5)node[midway,below] {\out{}};
            \draw [arrow] (n5) -- (n6);

            \begin{scope}[on background layer]
                \node[fit=(n3) (n4), rectangle, fill=red!30, draw=red, fill opacity=0.5, label={[label distance=-0.3cm,text depth=3ex]above:neural network}] (qm) {};
            \end{scope}
        \end{tikzpicture}
        }
    \caption{Regressor architecture.}
    \label{tik:regressor}
    \end{figure}
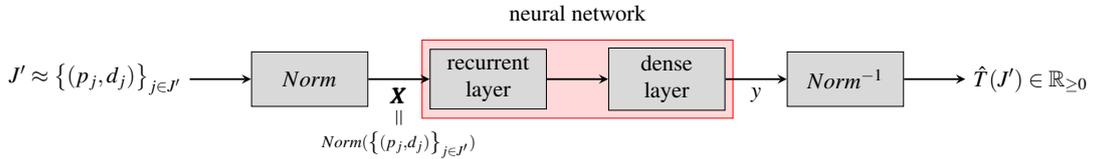

    The architecture of our regressor is illustrated in \autoref{tik:regressor}.
    Its input is job set $\jobsset'$ characterized by processing times and due dates of jobs. 
    Before the input is passed to the neural network, it is normalized to a matrix of features \inp{} by the block denoted $Norm$.
    The dimension of \inp{} is $2 \times |\jobsset'|$ where the first dimension corresponds to the two parameters of jobs $(p_j, d_j)$ entering the regressor. 
    The output of the neural network, $y$, is normalized as well; thus, it needs to be denormalized by $Norm^{-1}$ in order to get estimate $\hat{T}(\jobsset')$.

    The detailed description of the regressor is split into two subsections. 
    First, in \autoref{sec:reg:preprocessing} we  describe the normalization and denormalization of data, while the neural network is introduced in \autoref{sec:reg:nn}.

    \subsubsection{Normalization of the Input Data}
        \label{sec:reg:preprocessing}
        The normalization of the input data is a preprocessing step that improves the accuracy of the neural network, improves its numerical stability, and reduces training time.
        In our case, the normalization takes the job set $\jobsset^\prime$ and normalizes it to a  form suitable to the neural network.
        This preprocessing consists of (i) supplying the input to the network in a canonical form, (ii) normalization of the input features to $[0,1]$, and (iii) the normalization of the criterion.
        
        The normalization used in the regressor sorts the jobs in a defined order and scales the processing times, due dates, and the value of the objective function by suitable constants.
        It is important to note that the normalization of inputs for the neural network occurs in two different ways. 
        The first one concerns network training, when we need \inp{} and $y$, i.e., normalized input and output. 
        The other way occurs during the exploitation of the regressor, i.e., \emph{inference} from the neural network. 
        Then the input is normalized to \inp{} but the output of the network $y$ needs to be converted back to $\hat{T}(\jobsset')$ by block $Norm^{-1}$.
        The text below describes two normalization procedures that improved the precision of the regressor the most.
        
        Both normalizations sort $\jobsset^\prime$ in \edd order. The first normalization, denoted \normsumproc, scales down all $\proctime{j}, \duedate{j} :\forall j \in \jobsset^\prime$ as well as the objective function $\hat{T}(\jobsset')$ by factor $\max\left\{\sum_{\job \in \jobsset^\prime}\proctime{\job}, \max_{\job \in \jobsset^\prime}\duedate{\job} \right\}$.
        This normalizes the processing times and due dates to the interval $[0,1]$; nevertheless, the value of the objective function obtained for the rescaled parameters can be greater than 1. 
        Despite its simplicity, we have observed a noticeable impact of this normalization on the accuracy of the regressor.
        
        The second normalization, denoted as \gapeddinv, redefines the output of the neural network.
        The neural network is not trained to predict the target criterion value but rather the difference between the optimal objective value and the objective value achieved by \edd job ordering.
        This is the key innovation, as the neural network has an easier job modeling residue between the optimal and \edd solutions rather than modeling the objective function from scratch.
        Except for a different quantity to predict, \gapeddinv uses the same normalization constant for the processing times and due dates as \normsumproc; however, it redefines how the predicted value $\hat{T}(\jobsset')$ is computed.
        
        During the training phase of the network, the \gapeddinv normalization constructs the \edd sequence of jobs.
        Then it takes the optimal solution \optsequence{} and computes the optimality gap of the \edd sequence \eddsequence{} as $gap_{EDD} = \gapdeffinition{\eddsequence{}}{^\prime}$.
        Subsequently, to have the output of the neural network distributed in interval $[0,1]$, the normalized criterion is computed as $y = \frac{1}{1+gap_{EDD}}$.
        The inverse transformation, needed for the inference from the neural network during the \DHS{} run, proceeds the other way around. 
        First, it computes $gap_{EDD}=\frac{1}{y}-1$ from the normalized output $y$.
        Then it computes the total tardiness of the $\eddsequence{}$ sequence, which can be constructed in polynomial time. 
        Finally, it uses the definition of the \edd optimality gap to derive $\hat{T}(\jobsset')$.

    \subsubsection{Neural Network}
        \label{sec:reg:nn}
        The main reasons why we have used a deep neural network architecture in this paper are the following.
One of the main limitations of the application of classical machine learning models (including standard neural network) is that they expect a fixed-sized input, meaning that the dimension of the input data is always the same number, which needs to be chosen before the training and cannot be changed afterward.
However, the scheduling problem which we solve assumes an arbitrary number of jobs $\jobsset$, and it is not apparent how to compress the input instance into a fixed-sized input.
Note that this is radically different from, e.g., computer vision problems, where the input image can be naturally scaled down to match the target dimensions.

Another reason is that deep neural network architectures have sufficient capability to discover their own features of the input data.
For classical ML models, it is needed to design custom descriptors of the data (i.e., features).
In the case of our scheduling problem with total tardiness minimization, such features might be statistical indicators of processing times and dues (e.g., means and variances).
However, since we face an \nphard{} problem (i.e., even the prediction of the optimal objective is a hard problem), it is very unlikely that such simple descriptors would work well.
More complex interactions between the individual job parameters need to be considered for precise predictions.
Indeed, when we tested in our experiments simpler networks, we could see that the prediction accuracy is affected by the capacity (i.e., model complexity) of the prediction models. 
For example, in \cref{sec:exp:nn-params} we evaluate the effect of the complexity of the network on the quality of solutions.
There, we can see that reducing the size of the hidden state of the neural network has already a negative impact on the performance of \DHS{}.
However, even the features developed by the network with a limited complexity are still much more complex than the handcrafted features like the statistical ones mentioned above.
Thus, we let the network discover its own features from raw data.

        The solution to both these challenges is offered by recurrent neural networks~\cite{Sundermeyer2012}.
        A recurrent neural network is a type of network that can use previous outputs as inputs in successive iterations.
        It maintains a hidden state which is updated by the current input and the past value of the hidden state.
        The whole input to the recurrent neural network is represented by a sequence that is processed in an iterative fashion.
        The final output of the network is a function of the last hidden state.
        In our case, the input sequence is the instance of the scheduling problem presented in a predefined canonical order.
        In each step of the computation, the features (i.e., processing time and due date) of a single job are fed to the network to update its hidden state.
        Finally, when the sequence is processed, the hidden state of the network reflects the whole problem instance which can be used to make predictions, e.g., about its optimal objective value.
        
        The advantage of this model is the possibility of processing the input of any length (i.e., a variable number of jobs) and that the historical information is used during the whole computation over the sequence.
        This scheme has been shown to be effective for processing sequences with variable lengths, e.g., in \NLP.


        The neural network used inside the regressor consists of two layers (see the red box in \autoref{tik:regressor}).
        The first layer is a  recurrent neural network, which receives normalized job set \inp{} as the input.
        This layer is realized using the \LSTM architecture~\cit{Hochreiter1997}.
        The output of the recurrent layer, the hidden state, is passed into a \emph{dense layer} without an activation function, and it produces estimation \out{}.
        The main parameters of the recurrent layer are the size of the hidden state, also denoted as \emph{capacity} in the literature.
        The recurrent layer's capacity affects the amount of information which is the recurrent layer able to approximate.

        In our experiments, we compare two different types of recurrent layers, \LSTM and \GRU\cite{cho2014}. In general, \GRU is better suited for smaller training data sets. \LSTM is more general and has a more complex structure compared to \GRU.
        For example, \GRU has only one reset gate, which substitutes the function of update and the forget gate in \LSTM.
        The experimental results documented in Section~\ref{sec:experiments} show that \LSTM provides better results in our case.

        Training of the neural network is complicated by its integration into the decompositions.
        It means that the training and validation errors (used in the training phase of the neural network) are not computed in the same way as the testing error (measured during benchmarking of the entire \DHS{}).
        The training and validation errors of the neural network are measured in terms of the mean square error of predicted objective values on the training and validation samples.
        On the other hand, the testing error is measured as the optimality gap of \DHS{}, which exploits the neural network.
        Computation of training and validation errors reflecting the optimality gap obtained by the whole \DHS{} would be highly time demanding for the training phase.
        Therefore, our approach relies on the relation between the prediction of objective values and their use in Equation~\eqref{eq:exact-criterion} controlling decisions of \DHS{}.
        Thus the training is faster, and \DHS{} still provides accurate solutions.

    \subsection{Training Data Set Generation}    
        \label{sec:sol:data-gen}

        The training data set for a neural network usually consists of thousands of millions of training samples.
        However, producing a training data set with this size can be extremely time demanding in case of \nphard{} problems.
        Thus, it is extremely important to devise efficient methods that can produce it in a reasonable time.
        An equally important property of the training data set is its composition, i.e., how well it covers different parts of the parameter space.
        Essentially, the aim is to ensure that the distribution of training samples corresponds to the distribution of the test samples that arise from the decomposition during the run of \DHS{}.
        As it will be shown in the next section, the composition of samples in the training data set is indeed critical to the quality of solutions produced by \DHS{}.
        In this section, we propose two approaches to the generation of the training data set, and we describe their properties.
        
        To generate \SMTTP instances, we utilize the standard benchmark generator described by \authoretal{Potts} \cite{potts1991}.
        They generate processing times from a uniform distribution $[1,\maxproc]$ where \maxproc is the maximum processing time.
        The distribution of due dates is given by two parameters \rdd (range of due dates) and \tf (tardiness factor), describing the relations between the sum of processing times and the due dates of jobs.
        Specifically, each due date is drawn from a uniform distribution on interval $[(1-\tf-\rdd/2) \sum_{j \in J} p_j, (1-\tf+\rdd/2) \sum_{j \in J} p_j]$.
        This procedure is used in the literature to generate a benchmark set for measuring the quality of algorithms for \SMTTP and various related problems.
        However, we need to design an efficient procedure that creates the training data set with the following properties: (i) the training data set needs to contain training labels (i.e., optimal objective values) for all the samples, (ii) the distribution of training instances should reflect what will be encountered during the inference, and (iii) the procedure should be able to produce millions of training samples in a reasonable time.
        In the following lines, we describe two choices of such training data set generators with their properties.
        
        The most straightforward method is  \gensolve{}.
        It produces the training data set as follows.
        First, it generates a random instance by the generator described by \authoretal{Potts}~\cite{potts1991}, and then it is solved by \TTBR{} algorithm, which acquires the training label (optimal objective value) for that instance.
        The pair consisting of the instance and its solution is then treated as a single training sample.
        The entire training data set is obtained by generating random instances having a uniform distribution of $\numjobs$, $\rdd$, and $\tf$.
        
        {
        \updatemb{}
        \begin{table}[thb!]
         \caption{Time consumed for generation of a data set with $1.5 \cdot 10^6$ training samples.}
        \label{tab:exp:data-gen:generation-times}
        \centering
        \begin{subtable}[t]{0.25\textwidth}
        \caption{\gensolve{}}
        \centering
        \begin{tabular}{c|c}
                \toprule
                \n  & time [s] \\                
                \midrule
                $5$-$100$ & $16470$ \\
                $5$-$125$ & $17743$ \\
                $5$-$150$ & $20006$ \\
                $5$-$175$ & $23766$ \\
                $5$-$200$ & $29757$ \\
                \bottomrule
            \end{tabular}
        \end{subtable}
        \begin{subtable}[t]{0.25\textwidth}
        \caption{\gensubinstance{}}
        \centering
        \begin{tabular}{c|c}
                \toprule
                \n  & time [s] \\                
                \midrule
                $75$-$100$  & $594$ \\
                $100$-$125$ & $803$ \\
                $125$-$150$ & $832$ \\
                $150$-$175$ & $1040$ \\
                $175$-$200$ & $1020$ \\
                \bottomrule
            \end{tabular}
            
        \end{subtable}

        \end{table}
        }

        A different method to generate the training data set is named \gensubinstance.
        Its idea is to exploit all subproblems that are being solved by a decomposition-based algorithm (see Section~\ref{sec:smttpdec}) during the solution of a single problem instance.
        For the given input instance \jobssetrec and eligible position \position of the splitting job \jsplit both \edd, \spt decompositions generate two subproblems, i.e. \precprob{\jobssetrec,\position} and \follprob{\jobssetrec,\position}.
        Since the decomposition is applied recursively, these subproblems generate other subproblems.
        Thus, from a single run of a decomposition-based algorithm, multiple training instances together with their optimal objective values emerge.
        
        \gensubinstance proceeds as follows.
        First, it generates a problem instance by the generator described above.
        The instances are generated with constant \rdd and \tf (specifically $\rdd = 0.2$ and $\tf = 0.6$).
        The reason for this setting is that the instances with those parameters are the most computationally demanding compared to other $(\rdd, \tf)$~\cite{Shang2017}. 
        The second reason why it is possible to assume only a single $(\rdd, \tf)$ pair is that the decompositions generate problems \precprob{\jobssetrec,\position} and \follprob{\jobssetrec,\position} with a different $(\rdd, \tf)$ than the ones of \jobssetrec.
        Therefore, the generated data set covers different $(\rdd, \tf)$ pairs as well.  

        After an instance is generated, it is solved by a combination of \edd and \spt decomposition where the algorithm always selects the decomposition having the smallest $\overline{K^{\circ}}(\jobsset)$, $\circ\in\{\edd, \spt\}$.
        In each step, the instance is recursively decomposed into a set of subproblems \precprob{\jobssetrec,\position} and  \follprob{\jobssetrec,\position} for different eligible positions \position of the splitting job.
        Subsequently, all subproblems are solved by the recursive application of the decomposition.
        In the end, the optimal solutions of all the subproblems are known; thus, all pairs of subproblems and their solutions are put into the training data set.
        Therefore, one can get multiple training samples from a single problem instance solved to the optimality.
        In addition, newly generated subproblems differ in the number of jobs and their characteristics in terms of \textit{rdd} and \textit{tf} parameters, which enriches the composition of the training data set.
        
        Now, let us discuss the properties of the above generators.
        The \gensolve method has two main disadvantages.
        The first one is the time complexity of the generation.
        Every problem instance results in only a single sample of the training data set.
        On the other hand, the \gensubinstance naturally generates many training samples from a single input instance, which significantly reduces the time needed for the training data set synthesis.
        \cref{tab:exp:data-gen:generation-times} shows the times needed for generating data sets with $1.5\cdot 10^6$ training samples for different intervals of instance sizes.
        Comparing the two methods, it can be seen that \gensubinstance can generate the training data set more than twenty times faster.
        
        The second disadvantage of \gensolve method is the spectrum of generated instances. 
        While \gensolve method produces a uniform distribution of instances with respect to \n, \rdd, and \tf, \gensubinstance generates a spectrum of instances focused on the needs of \DHS{}.
        Indeed, the distribution of instances in the training data set produced by \gensubinstance is similar to the distribution of subproblems whose objective value needs to be estimated during the run of \DHS{} (which is far from being uniform in terms of \n, \rdd, and \tf).
        
         \begin{figure}[ht]
            \centering
           \begin{tikzpicture}
    \begin{axis}[
        width=0.8\textwidth,
        height=0.35\textwidth,
        grid=both,
        xlabel={$n$},
        ylabel={$\#$ training samples [-]},
        xmin=-1,
        xmax=200,
        scaled ticks = false,
        tick label style={
            /pgf/number format/fixed,
            /pgf/number format/precision=2
        },
        legend cell align={left},
        minor tick num=4,
        minor grid style={black!15},
        major grid style={black!30},
        label style={font=\footnotesize},
        tick label style={font=\scriptsize},
        legend pos=north east,
    ]

        \addplot[blue, thick, sharp plot] table [col sep=space] {distribution_subproblem_75.csv};
        \addplot[red, thick, sharp plot] table [col sep=space] {distribution_subproblem_125.csv};
        \addplot[green, thick, sharp plot] table [col sep=space] {distribution_subproblem_175.csv};
        \addplot[black, thick, sharp plot] table [col sep=space] {distribution_uniform.csv};
       \legend{\gensubinstance 75--100, \gensubinstance 125--150, \gensubinstance 175--200, \gensolve 5--200}
    \end{axis}
\end{tikzpicture}
            \caption{Distribution of training sample size for \gensubinstance method with different range of instances and \gensolve.}
            \label{fig:exp:data-gen:training-sample-size-distribution}
        \end{figure}
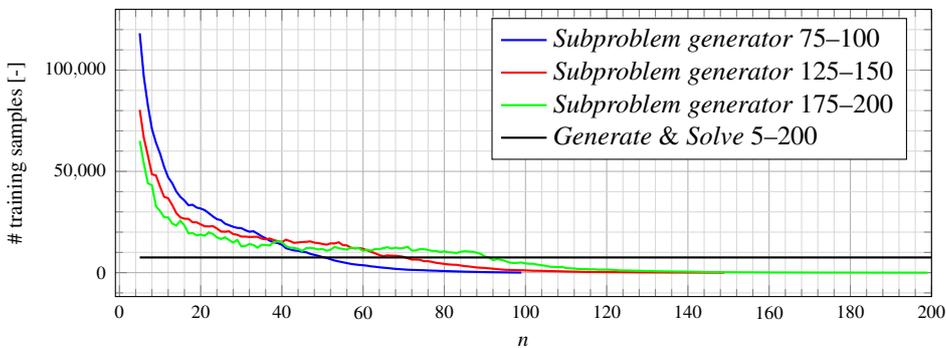


        The composition of instances in the training data set needs to be analyzed from two perspectives.
        This first one is the number of instances with a particular \n value. 
        The distribution of the training samples' sizes is shown in \cref{fig:exp:data-gen:training-sample-size-distribution} on a data set with $1.5\cdot 10^6$ training samples.
        The figure compares \gensolve for $n \in [5,200]$ and 
        \gensubinstance for three different ranges on $n$, i.e., $[75,100]$, $[125,150]$ and $[175,200]$.
        It shows that while the number of instances for different \n is constant for \gensolve, \gensubinstance generates significantly more instances with lower $n$. It indicates that decompositions inside \DHS{} need more examples with lower \n; thus, these instances are more important for the training of the neural network.


             \begin{figure}[thb!]
            \centering
           
                
            \begin{subfigure}{.51\textwidth}
                \centering
                \resizebox{\linewidth}{!}{
                    \input{img/dataset_density_graph/exp_dataset_size_n150-150x2O_rt_1608279941_984857_n150to151_rtm_csv_n0100}}%
                \vspace{-2mm}
                \caption{$\n = 100$}
            \end{subfigure}
            \hspace{-0.04\textwidth}
            \begin{subfigure}{.51\textwidth}
                \centering
                \resizebox{\linewidth}{!}{
                    \input{img/dataset_density_graph/exp_dataset_size_n150-150x2O_rt_1608279941_984857_n150to151_rtm_csv_n0060}}%
                \vspace{-2mm}
                \caption{$\n = 60$}
            \end{subfigure}
            \hspace{-0.04\textwidth}
            \begin{subfigure}{.51\textwidth}
                \centering
                \resizebox{\linewidth}{!}{
                    \input{img/dataset_density_graph/exp_dataset_size_n150-150x2O_rt_1608279941_984857_n150to151_rtm_csv_n0040}}%
                \vspace{-2mm}
                \caption{$\n = 40$}    
            \end{subfigure}%
             \hspace{-0.04\textwidth}
             \begin{subfigure}{.51\textwidth}
                \centering
                \resizebox{\linewidth}{!}{
                    \input{img/dataset_density_graph/exp_dataset_size_n150-150x2O_rt_1608279941_984857_n150to151_rtm_csv_nall}}%
                \vspace{-2mm}
                \caption{Aggregated over $\n \in\{5,150\}$}
            \end{subfigure}
            
            \caption[short]{Distribution of \rdd and \tf over \n in the data set generated by \gensubinstance.}
            \label{fig:sol:data-gen:rdd-tf-distribution}
        \end{figure}

        The second aspect that needs to be studied is the distribution of \rdd and \tf parameters of instances generated by both methods. 
        In this case, the situation is similar to the case with \n. 
        While instances with defined $(\rdd, \tf)$ parameters are uniformly distributed in the training data set produced by \gensolve method, the distribution of $(\rdd, \tf)$ is much more uneven in the case of \gensubinstance which is demonstrated in \cref{fig:sol:data-gen:rdd-tf-distribution}.
        It illustrates the frequency of samples with a particular value of $(\rdd, \tf)$ arising from subproblems generated from 20 instances having $150$ jobs with $\rdd = 0.2$ and $\tf = 0.6$.
        The samples are divided into particular categories according to the number of jobs \n.
        
        For $\n = 100$, one can see that the instances are close to the initial setting of $\rdd = 0.2$ and $\tf = 0.6$.
        With the decreasing value of \n (which corresponds to smaller subproblems), the distribution of \rdd and \tf shifts significantly. 
        The mass of the samples for the whole data set drifts towards $\rdd = 0.4$ and $\tf = 0.5$ and the covariance of $(\rdd, \tf)$ increases as well.
        These observations underline that the composition of instances that occur during a single run of \DHS{} is significantly different from a uniform distribution produced by \gensolve.
        At the same time, it shows clear advantages of data sets produced by \gensubinstance since their distribution is closer to what \DHS{} encounters.
        For more details, see the comparison of these two approaches in \cref{tab:exp:sota:p5000}, \cref{fig:exp:best-p100}, and \cref{fig:exp:best-p5000} in \cref{sec:exp:comparison}.



    \subsection{Time Complexity of the Scheduling Algorithm}
        \label{sec:dec:time-complexity}
        In this section, we present an analysis of the worst-case run time of \DHS{}.    
        The most time-consuming part of \DHS{} is the estimation of the optimal position for the splitting job \longestjobposition.
        To do so, the algorithm tests at most \n positions. For each position the algorithm uses the regressor to compute estimates \precprob{\jobssetrec,\position} and \follprob{\jobssetrec,\position} having $O(\n)$ time complexity.
        Thus, the estimation of the optimal position takes $O(\n^2)$.
        In the worst case, when the decomposition removes only one job from \jobsset in every recursion, \DHS{} makes $O(\n)$ selections of position \longestjobposition.
        Therefore, the worst-case time complexity of \DHS{} is $O(n^3)$.
        Nevertheless, as the experiments in \cref{sec:exp:comparison} show, the algorithm is very fast.
        The average CPU time on instances with 800 jobs is below 15 s.

\section{Experimental Results}
\label{sec:experiments}
In this section, we present the experimental results related to \DHS{} focused on its relation to the neural network used in the regressor. 
The section is organized as follows.
First, we describe the hardware and software used for the evaluation and the method of generating testing instances.
Next, we present the comparison of our approach with the \SOTA heuristic, the \SOTA exact algorithm, and results from our previous work \cite{bouska2019}.
The following sections describe a detailed benchmarking of the algorithm.
\cref{sec:sol:data-gen} presents the experiments with the generation of the training data set.
The influence of neural network hyperparameters on \DHS{} is investigated in \cref{sec:exp:nn-params}, and
\cref{sec:exp:horda-params} describes sensitivity of \DHS{} to its parameters.

\subsection{Experimental setup}
    \label{sec:exp:setup}
    Experiments were run on a single CPU core of the Xeon Gold 6140 processor with a memory limit set to 8 GB of RAM.
    \DHS{} and \NBR{} algorithms were implemented in Python 3.7, and the neural network is implemented using TensorFlow 2 trained on a GPU Nvidia GTX 1080 Ti.
    Source code of \DHS{} is available at \href{https://github.com/CTU-IIG/horda}{GitHub repository (https://github.com/CTU-IIG/horda)};
    Source code of \TTBR{} algorithm was provided by authors of \cit{Garraffa2018} and is implemented in C++.
    
    Instances used in this paper were generated in the manner suggested by Potts and Van Wassenhove~\cite{potts1991}.
    They generate processing times of jobs uniformly on the interval from $1$ to $100$.
    Since we want to study the impact of processing time on the quality of our algorithm, we define maximal processing time \maxproc and generate the processing time of jobs uniformly on the interval from $1$ to the \maxproc.
    The analysis presented in~\cite{Shang2017} implies that the hardest instances occur for $\rdd = 0.2$ and $\tf = 0.6$; therefore, our experiments focus on them.
    We use two data sets: the first with $\maxproc = 100$, and the second, with $\maxproc = 5000$, to mimic the duration of jobs in seconds used, e.g., nowadays in Advanced Planning and Scheduling systems.
    In all experiments, we have used separate data sets for training and evaluation; thus, all the methods are tested on instances that were not seen during the training.
    
    For all graphs presented in this section, the points in the figure represent the mean over all instances of the same size. The colored surface represents the standard deviation of the measurement, and the line represents the running mean of the last five values.
    The evaluation data sets consist of 20 randomly generated instances for each \n divisible by 5.
    The results are represented in a form of the optimality gap defined as $\gapdeffinition{\sequence}{}\cdot 100$, where \sequence{} is a heuristic solution and \objopt{\jobsset} is the optimal objective value.
    Next, note that the training of neural networks is a random process (e.g., data set shuffling, initial weights, numerical instability), and the results vary over different runs.
    Thus, for each settings of neural network hyperparameters, we trained all neural networks five times with the same parameters and used the best one to carry out the experiments.
    
\subsection{Comparison with \SOTA approaches}
    \label{sec:exp:comparison}

    In this section, we present a comparison of our results with \SOTA approaches in terms of the optimality gap and run time.
    Our results are compared with \NBR{} \cit{holsenback1992} \SOTA heuristic, its combination inside the decomposition \DHS{\NBR}{} \cite{della2004}, \TTBR{} \cit{Garraffa2018} as the \SOTA exact algorithm, and our previous \ML based approach denoted as \DHS{NN2020} \cite{bouska2019}.
    
    For a fair comparison of all methods, we limited the run time of \TTBR{} to 10\jed{s}, and 15\jed{s}. Time limit 15\jed{s} corresponds to the maximum time needed by \DHS{} to solve the largest instances.
    {\updatemb{}
    Notice that this gives an advantage to \TTBR{} on small instances, but it makes the comparison simpler.
    Thus the comparison is made in favor of \TTBR{}.
    We report these as \TTBR{10}, and \TTBR{15} respectively.
    }
    
    Results of the presented approach are also compared with our previous approach \DHS{NN2020}, presented in \cite{bouska2019}. It utilizes \DHS{} heuristic with the \LSTM based neural network (with capacity 256 and \normsumproc criterion normalization) trained on a data set generated by \gensolve (see \cref{sec:sol:data-gen}).

    \begin{table}[]
        \centering
        \caption{Comparison with \SOTA approaches on instances with $\maxproc = 5000.$}
        \label{tab:exp:sota:p5000}
        \small
        \begin{tabular}{l |r | c | c | c | c | c c | c | c c}
            \toprule
            $n$     & \makecell{\TTBR{}\\\textsc{\lowercase{($\infty$ s)}}}  & \makecell{\TTBR{}\\\textsc{\lowercase{(10 s)}}} & \makecell{\TTBR{}\\\textsc{\lowercase{(15 s)}}} & \textsc{nbr} & \makecell{\DHS{}\\\textsc{\lowercase{(MDD)}}} & \multicolumn{2}{c|}{\makecell{\DHS{}\\\textsc{(nbr)}}} &
            \makecell{\DHS{}\\\textsc{\lowercase{(NN2020)}}} & \multicolumn{2}{c}{\makecell{\DHS{}\\\textsc{(best)}}}\\
             &\scriptsize time [s] &\scriptsize gap [\%] &\scriptsize gap [\%] &\scriptsize gap [\%]  &\scriptsize gap [\%]  &\scriptsize gap [\%] &\scriptsize time [s] &\scriptsize gap [\%]     &\scriptsize gap [\%] &\scriptsize time [s] \\
            \midrule
            
            \multirow{2}{*}{\scriptsize $5-45$ } &$0.02$ & $0.00$ & $\textbf{0.00}$ & $0.95$ & $0.38$ & $0.09$ & $0.01$ & $0.56$ & $0.22$ & $0.03$ \\ [-6pt]
 &\scriptsize $\pm0.00$ &\scriptsize $\pm0.00$ &\scriptsize $\pm0.00$ &\scriptsize $\pm2.92$ &\scriptsize $\pm0.87$ &\scriptsize $\pm0.28$ &\scriptsize $\pm0.01$ &\scriptsize $\pm0.86$ &\scriptsize $\pm0.51$ &\scriptsize $\pm0.10$ \\
\multirow{2}{*}{\scriptsize $50-95$ } &$0.03$ & $0.00$ & $\textbf{0.00}$ & $1.22$ & $0.97$ & $0.36$ & $0.09$ & $0.64$ & $0.22$ & $0.18$ \\ [-6pt]
 &\scriptsize $\pm0.03$ &\scriptsize $\pm0.00$ &\scriptsize $\pm0.00$ &\scriptsize $\pm0.91$ &\scriptsize $\pm1.38$ &\scriptsize $\pm0.48$ &\scriptsize $\pm0.05$ &\scriptsize $\pm0.42$ &\scriptsize $\pm0.29$ &\scriptsize $\pm0.08$ \\
\multirow{2}{*}{\scriptsize $100-145$ } &$0.34$ & $0.00$ & $\textbf{0.00}$ & $1.59$ & $1.11$ & $0.71$ & $0.31$ & $0.46$ & $0.39$ & $0.46$ \\ [-6pt]
 &\scriptsize $\pm0.55$ &\scriptsize $\pm0.00$ &\scriptsize $\pm0.00$ &\scriptsize $\pm0.70$ &\scriptsize $\pm1.04$ &\scriptsize $\pm0.48$ &\scriptsize $\pm0.13$ &\scriptsize $\pm0.31$ &\scriptsize $\pm0.40$ &\scriptsize $\pm0.16$ \\
\multirow{2}{*}{\scriptsize $150-195$ } &$2.30$ & $0.01$ & $\textbf{0.00}$ & $1.86$ & $1.39$ & $0.92$ & $0.66$ & $0.47$ & $0.50$ & $0.93$ \\ [-6pt]
 &\scriptsize $\pm2.34$ &\scriptsize $\pm0.07$ &\scriptsize $\pm0.03$ &\scriptsize $\pm0.64$ &\scriptsize $\pm1.20$ &\scriptsize $\pm0.45$ &\scriptsize $\pm0.22$ &\scriptsize $\pm0.33$ &\scriptsize $\pm0.45$ &\scriptsize $\pm0.22$ \\
\multirow{2}{*}{\scriptsize $200-245$ } &$11.09$ & $0.17$ & $\textbf{0.07}$ & $1.99$ & $1.54$ & $1.05$ & $1.25$ & $0.59$ & $0.45$ & $1.56$ \\ [-6pt]
 &\scriptsize $\pm10.69$ &\scriptsize $\pm0.30$ &\scriptsize $\pm0.18$ &\scriptsize $\pm0.65$ &\scriptsize $\pm1.12$ &\scriptsize $\pm0.48$ &\scriptsize $\pm0.41$ &\scriptsize $\pm0.29$ &\scriptsize $\pm0.29$ &\scriptsize $\pm0.31$ \\
\multirow{2}{*}{\scriptsize $250-295$ } &$38.37$ & $0.82$ & $0.56$ & $2.08$ & $1.60$ & $1.19$ & $2.07$ & $0.55$ & $\textbf{0.37}$ & $2.30$ \\ [-6pt]
 &\scriptsize $\pm27.57$ &\scriptsize $\pm0.68$ &\scriptsize $\pm0.57$ &\scriptsize $\pm0.52$ &\scriptsize $\pm1.13$ &\scriptsize $\pm0.44$ &\scriptsize $\pm0.66$ &\scriptsize $\pm0.27$ &\scriptsize $\pm0.24$ &\scriptsize $\pm0.39$ \\
\multirow{2}{*}{\scriptsize $300-345$ } &$93.19$ & $1.32$ & $1.02$ & $2.29$ & $1.75$ & $1.33$ & $2.96$ & $0.57$ & $\textbf{0.31}$ & $3.07$ \\ [-6pt]
 &\scriptsize $\pm59.73$ &\scriptsize $\pm0.80$ &\scriptsize $\pm0.72$ &\scriptsize $\pm0.54$ &\scriptsize $\pm1.19$ &\scriptsize $\pm0.45$ &\scriptsize $\pm0.88$ &\scriptsize $\pm0.35$ &\scriptsize $\pm0.17$ &\scriptsize $\pm0.43$ \\
\multirow{2}{*}{\scriptsize $350-395$ } &$209.70$ & $1.99$ & $1.68$ & $2.35$ & $2.06$ & $1.42$ & $4.41$ & $1.16$ & $\textbf{0.27}$ & $4.03$ \\ [-6pt]
 &\scriptsize $\pm113.55$ &\scriptsize $\pm0.87$ &\scriptsize $\pm0.81$ &\scriptsize $\pm0.48$ &\scriptsize $\pm1.34$ &\scriptsize $\pm0.45$ &\scriptsize $\pm1.25$ &\scriptsize $\pm0.60$ &\scriptsize $\pm0.16$ &\scriptsize $\pm0.51$ \\
\multirow{2}{*}{\scriptsize $400-445$ } &$464.90$ & $2.62$ & $2.33$ & $2.36$ & $2.31$ & $1.51$ & $6.00$ & $1.72$ & $\textbf{0.24}$ & $5.01$ \\ [-6pt]
 &\scriptsize $\pm246.80$ &\scriptsize $\pm0.87$ &\scriptsize $\pm0.83$ &\scriptsize $\pm0.45$ &\scriptsize $\pm1.34$ &\scriptsize $\pm0.41$ &\scriptsize $\pm1.75$ &\scriptsize $\pm0.62$ &\scriptsize $\pm0.15$ &\scriptsize $\pm0.63$ \\
\multirow{2}{*}{\scriptsize $450-495$ } &$814.21$ & $2.83$ & $2.55$ & $2.43$ & $2.33$ & $1.56$ & $7.70$ & $1.32$ & $\textbf{0.22}$ & $6.10$ \\ [-6pt]
 &\scriptsize $\pm398.34$ &\scriptsize $\pm0.83$ &\scriptsize $\pm0.78$ &\scriptsize $\pm0.45$ &\scriptsize $\pm1.15$ &\scriptsize $\pm0.41$ &\scriptsize $\pm2.15$ &\scriptsize $\pm0.60$ &\scriptsize $\pm0.13$ &\scriptsize $\pm0.79$ \\
\multirow{2}{*}{\scriptsize $500-545$ } &$1528.93$ & $3.23$ & $2.98$ & $2.43$ & $2.34$ & $1.54$ & $10.28$ & $1.28$ & $\textbf{0.20}$ & $7.19$ \\ [-6pt]
 &\scriptsize $\pm620.02$ &\scriptsize $\pm0.80$ &\scriptsize $\pm0.80$ &\scriptsize $\pm0.43$ &\scriptsize $\pm1.24$ &\scriptsize $\pm0.38$ &\scriptsize $\pm2.88$ &\scriptsize $\pm0.51$ &\scriptsize $\pm0.14$ &\scriptsize $\pm0.80$ \\
\multirow{2}{*}{\scriptsize $550-595$ } &$2578.16$ & $3.49$ & $3.27$ & $2.52$ & $2.34$ & $1.66$ & $12.73$ & $1.20$ & $\textbf{0.18}$ & $8.50$ \\ [-6pt]
 &\scriptsize $\pm1128.52$ &\scriptsize $\pm0.71$ &\scriptsize $\pm0.70$ &\scriptsize $\pm0.41$ &\scriptsize $\pm1.12$ &\scriptsize $\pm0.44$ &\scriptsize $\pm3.00$ &\scriptsize $\pm0.48$ &\scriptsize $\pm0.12$ &\scriptsize $\pm0.89$ \\
\multirow{2}{*}{\scriptsize $600-645$ } &$4436.89$ & $3.70$ & $3.46$ & $2.54$ & $2.37$ & $1.62$ & $16.13$ & $1.20$ & $\textbf{0.16}$ & $9.82$ \\ [-6pt]
 &\scriptsize $\pm2341.48$ &\scriptsize $\pm0.70$ &\scriptsize $\pm0.70$ &\scriptsize $\pm0.37$ &\scriptsize $\pm1.15$ &\scriptsize $\pm0.33$ &\scriptsize $\pm3.90$ &\scriptsize $\pm0.41$ &\scriptsize $\pm0.08$ &\scriptsize $\pm0.97$ \\
\multirow{2}{*}{\scriptsize $650-695$ } &$7311.81$ & $3.86$ & $3.65$ & $2.54$ & $2.53$ & $1.66$ & $19.94$ & $1.05$ & $\textbf{0.15}$ & $11.18$ \\ [-6pt]
 &\scriptsize $\pm3689.10$ &\scriptsize $\pm0.79$ &\scriptsize $\pm0.78$ &\scriptsize $\pm0.40$ &\scriptsize $\pm1.23$ &\scriptsize $\pm0.37$ &\scriptsize $\pm4.92$ &\scriptsize $\pm0.50$ &\scriptsize $\pm0.09$ &\scriptsize $\pm1.08$ \\
\multirow{2}{*}{\scriptsize $700-745$ } &$11390.34$ & $3.90$ & $3.71$ & $2.56$ & $2.57$ & $1.74$ & $22.96$ & $1.01$ & $\textbf{0.16}$ & $12.71$ \\ [-6pt]
 &\scriptsize $\pm4360.18$ &\scriptsize $\pm0.75$ &\scriptsize $\pm0.73$ &\scriptsize $\pm0.41$ &\scriptsize $\pm1.09$ &\scriptsize $\pm0.38$ &\scriptsize $\pm6.18$ &\scriptsize $\pm0.37$ &\scriptsize $\pm0.09$ &\scriptsize $\pm1.20$ \\
\multirow{2}{*}{\scriptsize $750-795$ } &$14135.70$ & $3.92$ & $3.75$ & $2.59$ & $2.64$ & $1.65$ & $28.79$ & $0.89$ & $\textbf{0.11}$ & $13.95$ \\ [-6pt]
 &\scriptsize $\pm4814.88$ &\scriptsize $\pm0.64$ &\scriptsize $\pm0.64$ &\scriptsize $\pm0.31$ &\scriptsize $\pm0.91$ &\scriptsize $\pm0.37$ &\scriptsize $\pm6.96$ &\scriptsize $\pm0.41$ &\scriptsize $\pm0.08$ &\scriptsize $\pm1.29$ \\
\midrule
\multirow{2}{*}{\scriptsize $avg$ } &$2688.50$ & $1.99$ & $1.81$ & $2.14$ & $1.89$ & $1.25$ & $8.52$ & $0.92$ & $\textbf{0.26}$ & $5.44$ \\ [-6pt]
 &\scriptsize $\pm1113.36$ &\scriptsize $\pm0.55$ &\scriptsize $\pm0.52$ &\scriptsize $\pm0.66$ &\scriptsize $\pm1.16$ &\scriptsize $\pm0.41$ &\scriptsize $\pm2.21$ &\scriptsize $\pm0.46$ &\scriptsize $\pm0.21$ &\scriptsize $\pm0.61$ \\

            \bottomrule
        \end{tabular}
    
    \end{table}

        \begin{table}[]
        \centering
        \caption{Comparison with a Genetic Algorithm~\cite{suer2012}.}
        \label{tab:exp:ga}
        \small
        \begin{tabular}{l |r | c c | c c}
            \toprule
               & \TTBR{}  & \multicolumn{2}{c
               |}{GA~\cite{suer2012}} & \multicolumn{2}{c}{\DHS{best}}\\
            $n$ &\scriptsize time [s] &\scriptsize gap [\%] &\scriptsize time [s] & \scriptsize gap [\%] &\scriptsize time [s] \\
            \midrule
            \multirow{2}{*}{\scriptsize $10$ }  & $0.011$ & $0$   & $2$   & $0$     & $0$  \\ [-6pt]
            & {\scriptsize $\pm 0.001$} & & &{\scriptsize $\pm 0$} &{\scriptsize $\pm 0$ }\\
            \multirow{2}{*}{\scriptsize $20$ }  & $0.012$ & $0$   & $51$  & $0$     & $0.06$\\  [-6pt]
            & {\scriptsize $\pm 0.002$} & & &{\scriptsize $\pm 0$} &{\scriptsize $\pm 0.004$ }\\
            \multirow{2}{*}{\scriptsize $30$ }  & $0.012$ & $0$   & $354$ & $0.101$ & $0.013$\\  [-6pt]
            & {\scriptsize $\pm 0.002$} & & &{\scriptsize $\pm 0.101$} &{\scriptsize $\pm 0.005$ }\\
            \multirow{2}{*}{\scriptsize $50$ }  & $0.013$ & $2.12$& $536$ & $0.006$ & $0.018$\\  [-6pt]
            & {\scriptsize $\pm 0.002$} & & &{\scriptsize $\pm 0.006$} &{\scriptsize $\pm 0.007$ }\\
            \multirow{2}{*}{\scriptsize $100$ } & $0.013$ & $6.32$& $1083$& $0.002$ & $0.044$\\  [-6pt]
             & {\scriptsize $\pm 0.002$} & & &{\scriptsize $\pm  0.002$} &{\scriptsize $\pm 0.014$ }\\
            \bottomrule
        \end{tabular}
    
    \end{table}

    \DHS{best} method, proposed in this paper, uses the regressor, as it is presented in \cref{tik:regressor}.
    The recurrent layer is \LSTM with a capacity of $256$.
    The neural network uses \gapeddinv normalization, and it is trained on a data set generated with \gensubinstance. The data set was generated from input instances having $\numjobs \in [75,100]$ while for each $\numjobs$ there were $20$ instances used by \gensubinstance to generate the training data set. In total, the data set consists of $1.6\cdot10^6$ instances.

    First, let us discuss the time needed to train the neural network used in \DHS{best}.
    The generation of the training data set takes about $600\jed{s}$, and it is discussed in detail in \cref{sec:sol:data-gen}.
    The training time of the neural network is about $3$ hours.
    It comes as no surprise that the neural network training time is large compared to the time needed to solve a single instance by \TTBR{}.
    The approach presented here assumes that the time needed for the neural network training is a part of the algorithm development, but the online execution of the algorithm must be fast.
    Indeed, the training time is spent just once, whereas the inference from the neural network after the training is much faster.

    \cref{tab:exp:sota:p5000} compares the run times and optimality gaps of \TTBR{}, \TTBR{10}, \TTBR{15}, \NBR, \DHS{mdd}, \DHS{nbr}, \DHS{NN2020}, and \DHS{best} on instances generated with $\maxproc = 5000$.
    The bold numbers in the table represent the best result on the interval of \n in terms of the optimality gap.
    The time required for computing the optimal solution by \TTBR{} is shown in the second column of the table.
    For the largest instances with $\maxproc = 5000$, the computation takes up to $15\powten{3}\jed{s}$.
    Due to a software issue, we were not able to solve larger instances by \TTBR{}.
    \TTBR{} with limited time (\TTBR{10}, \TTBR{15}, i.e., third and fourth column, respectively) has small gaps on small instances, but from $\n = 250$ it perform worse than \DHS{best}. The average gap on the biggest instance is $3.92\%$ for \TTBR{10}, and $3.75\%$ for \TTBR{15}.
    \NBR heuristic has the optimality gap almost independent of the size of instances achieving an average gap $2.14\%$ over all instances (see the fifth column in the table). 
    Nevertheless, this result is significantly worst compared to the gap obtained by \DHS{best}.
    On the other hand, the average run time of \NBR heuristic for the biggest instances is around $0.85\jed{s}$ compared to \DHS{best} achieving $13.95\jed{s}$.
    
    As noted by Della Croce \emph{et al.} \cite{della2004}, constructive heuristics such as \textsc{mdd} \cite{bertrand1982} and \NBR can be embedded into the decomposition utilized by \DHS{} as well. 
    Thus, \DHS{mdd} and \DHS{nbr} columns present the integration of the respective constructive heuristic into \DHS{}, where it acts as a regressor replacing the neural network.
    The embedded \textsc{mdd} heuristic has an average gap $1.89\%$ (over all instances), and the quality is slightly better than \NBR.
    The integration of \NBR into \DHS{nbr} improves the gap from $2.14\%$ to $1.25\%$ (over all instances).
    Nevertheless, \DHS{best} has a significantly lower run time and almost five times lower gap.
    The table illustrates that within the 15 seconds time limit for instances having more than 250 jobs, the best results were achieved by \DHS{best}.

    \cref{tab:exp:ga} compares the results of \TTBR{}, \DHS{best} and a genetic algorithm presented by \authoretal{Suer} \cite{suer2012}.
    The \authoretal{Suer} defines its own evaluation protocol to generate test instances.
    Therefore, we evaluate \TTBR{} and \DHS{best} on instances generated in the same manner as suggested in \cite{suer2012} and show the results in the separate table.
    The run time of the genetic algorithm is scaled by the power ratio between the processor used by \authoretal{Suer} and us, i.e., $0.81$.
    The run time of \TTBR{} is $0.013\jed{s}$ for instances with $100$ jobs, while the genetic algorithm has about $1083\jed{s}$ (see column GA in the table).
    Our method obtains the solutions with an average gap $0.002\%$ within the $0.044\jed{s}$ for instances with $100$ jobs, which is superior to the results reported by \authoretal{Suer}

    \begin{figure*}[htb!]
        \centering
        \captionsetup{skip=0pt}
        \input{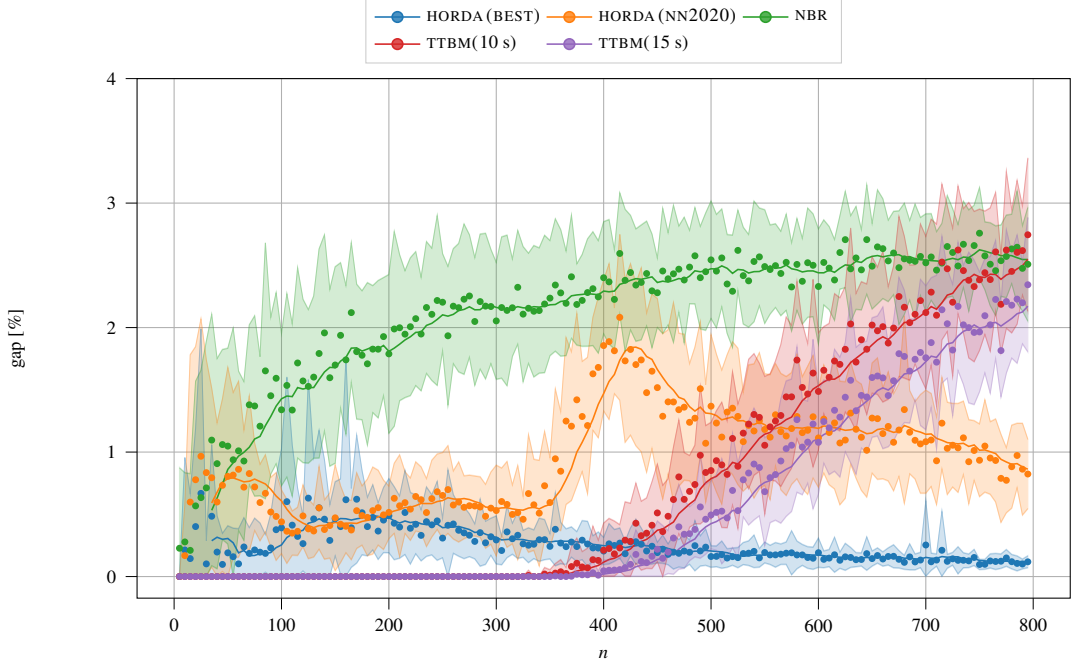}
        \caption{Optimality gap for instances with $\maxproc = 100.$}
        \label{fig:exp:best-p100}
    \end{figure*}
    \begin{figure*}[htb!]
        \centering
        \captionsetup{skip=0pt}
        \input{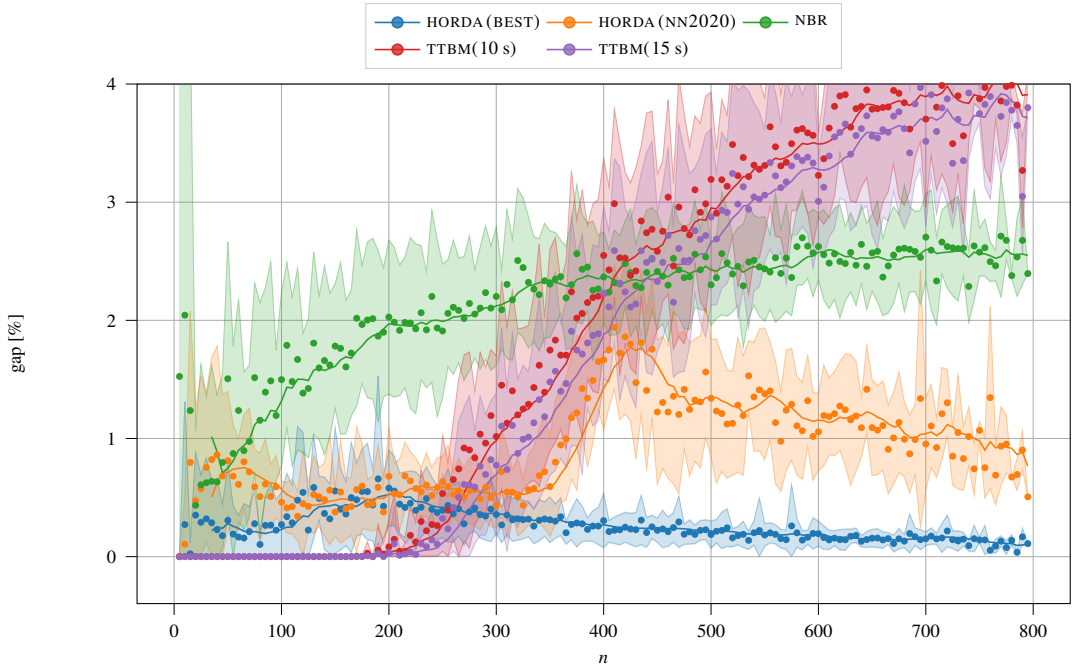}
        \caption{Optimality gap for instances with $\maxproc = 5000.$}
        \label{fig:exp:best-p5000}
    \end{figure*}

    The dependency of the optimality gap on the number of jobs for $\maxproc = 100$ and $\maxproc = 5000$ is shown in \cref{fig:exp:best-p100} and \cref{fig:exp:best-p5000}, respectively.
    The figures show the optimality gap of \NBR{}, \TTBR{10}, \TTBR{15}, \DHS{NN2020}, and \DHS{best} on \numjobs depending on \numjobs from $5$ to $800$.
    By comparing these two graphs, one can see that results of \TTBR{} are the best on the smaller instances.
    On the other hand, the \TTBR{} run time is heavily dependent on \maxproc; while other methods do not depend on \maxproc, which provides an advantage to them.
    Our previous method \DHS{NN2020} has very good results for instances with no more than $350$ jobs, and its gap is less than $1\%$.
    The method presented in this paper, i.e., \DHS{best} has the optimality gap for all sizes of instances under $0.5\%$, i.e., it is superior to all other methods for instances with about more than $450$ jobs for $\maxproc=100$, and $250$ jobs for $\maxproc=5000$.
    Neither the run time nor the gap of \DHS{best} depends on the \maxproc.
    Indeed, thanks to the normalization of the input data and the generalization capabilities of the neural network, it is possible to train it on the data set with $\maxproc = 100$ and apply it on instances generated with $\maxproc = 5000$.
    Thus, even thought \DHS{best} is trained on instances generated for $\maxproc = 100$, the results in terms of optimality gap are almost the same for the $\maxproc = 5000$ as for $\maxproc = 100$.
    This provides us some confidence that \DHS{best} is able to generalize for instances with different \maxproc.
    Moreover, an advantage of our approach is that the training data set can be generated for $\maxproc = 100$, which is much faster than the generation of the training data set with $\maxproc = 5000$.
    


\subsection{Neural network hyperparameters}
\label{sec:exp:nn-params}
    \begin{figure}[htb!]
        \centering
        \captionsetup{skip=0pt}
        \input{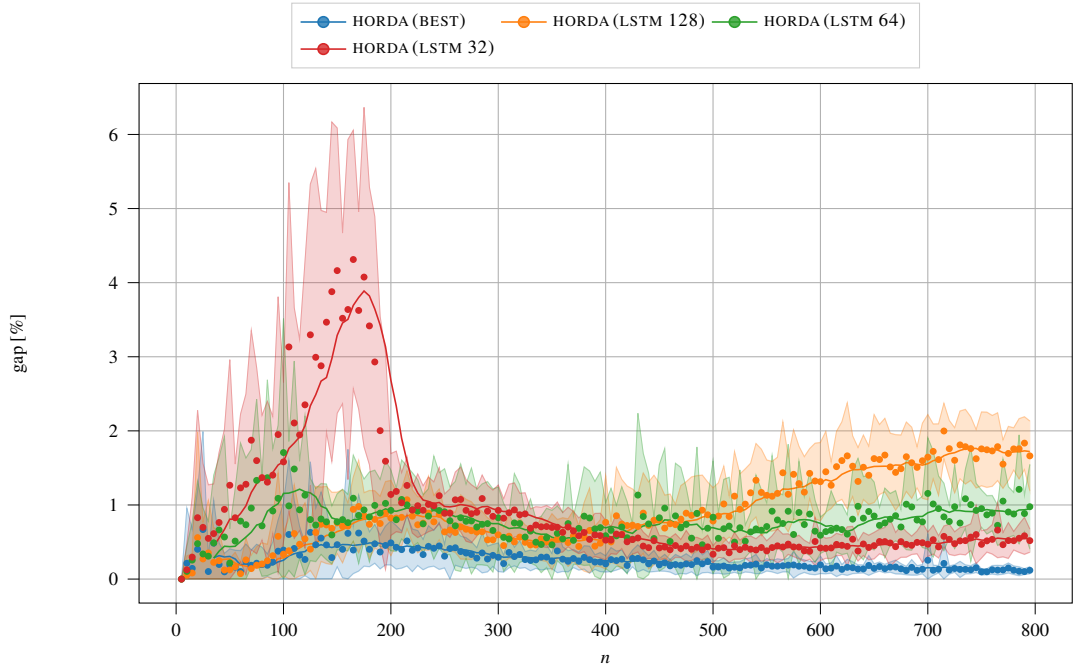}
        \caption{Optimality gap of \DHS{} with different capacity of \LSTM as a regressor.}
        \label{fig:exp:nn:capacity-gap}
    \end{figure}
    \begin{figure}[htb!]
        \centering
        \captionsetup{skip=0pt}
        \input{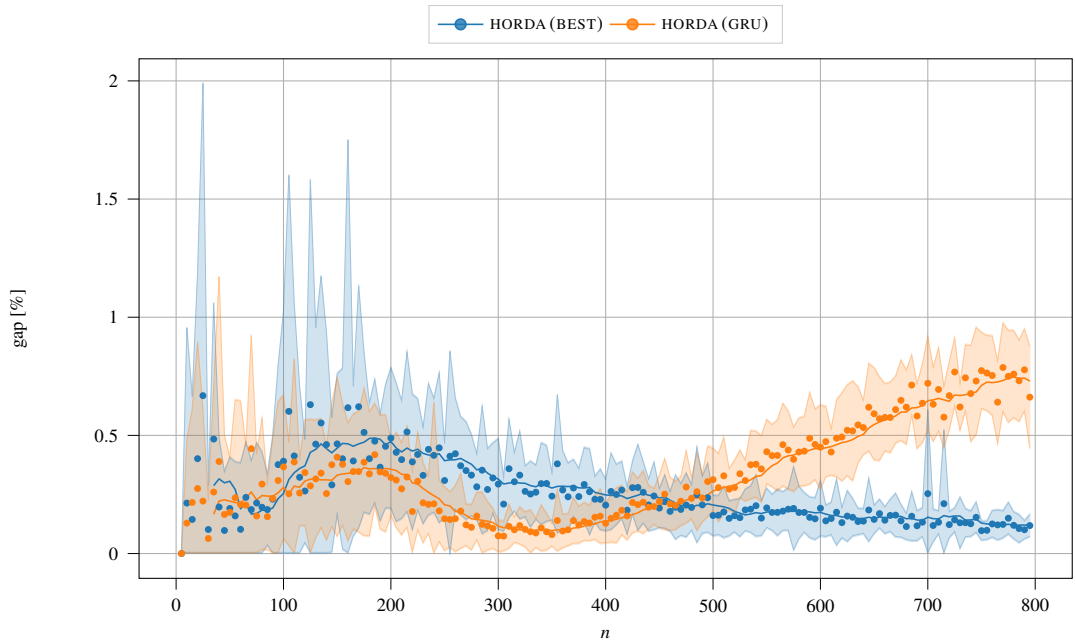}
        \caption{Optimality gap of \DHS{} with different neural networks.}
        \label{fig:exp:nn:capacity-gru}
    \end{figure}
    \begin{figure}[htb!]
        \centering
        \captionsetup{skip=0pt}
        \input{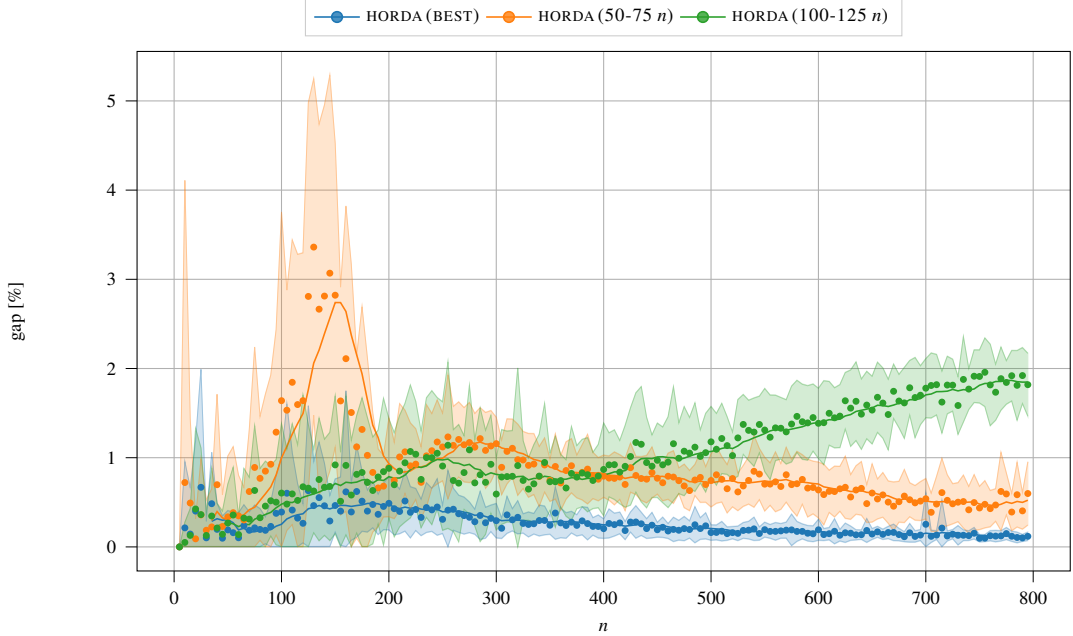}
        \caption{Optimality gap of \DHS{} with  a different number of jobs in training instances.}
        \label{fig:exp:nn:data-size:move-n}
    \end{figure}
    \begin{figure}[htb!]
        \centering
        \captionsetup{skip=0pt}
        \input{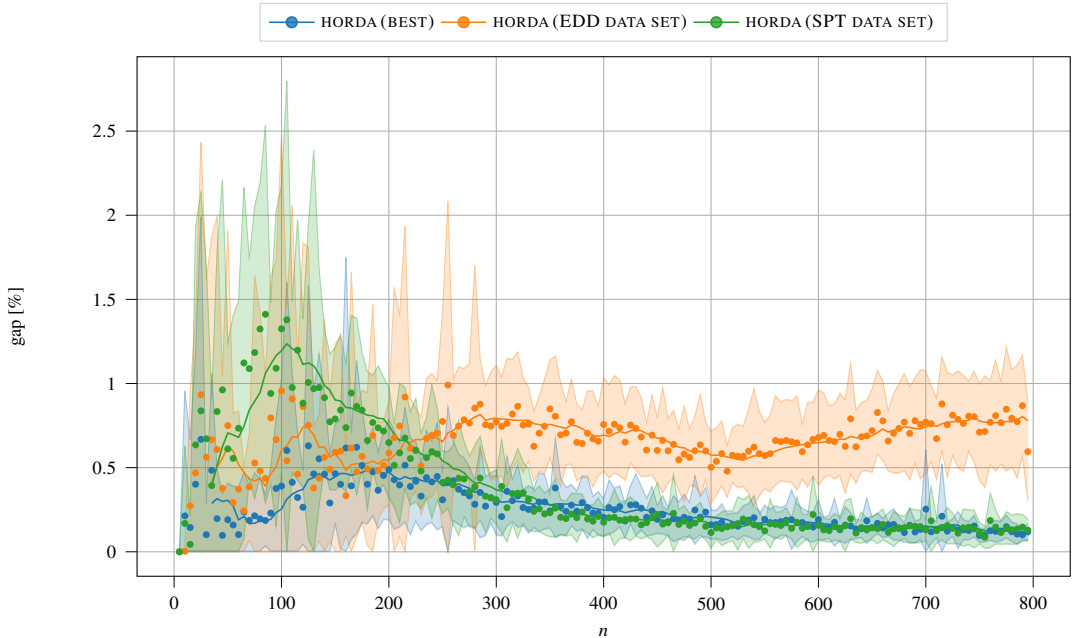}
        \caption{Optimality gap of \DHS{} with different decomposition used to generating training data set.}
        \label{fig:exp:nn:data-gen}
    \end{figure}
    \begin{figure}[htb!]
        \centering
        \captionsetup{skip=0pt}
        \input{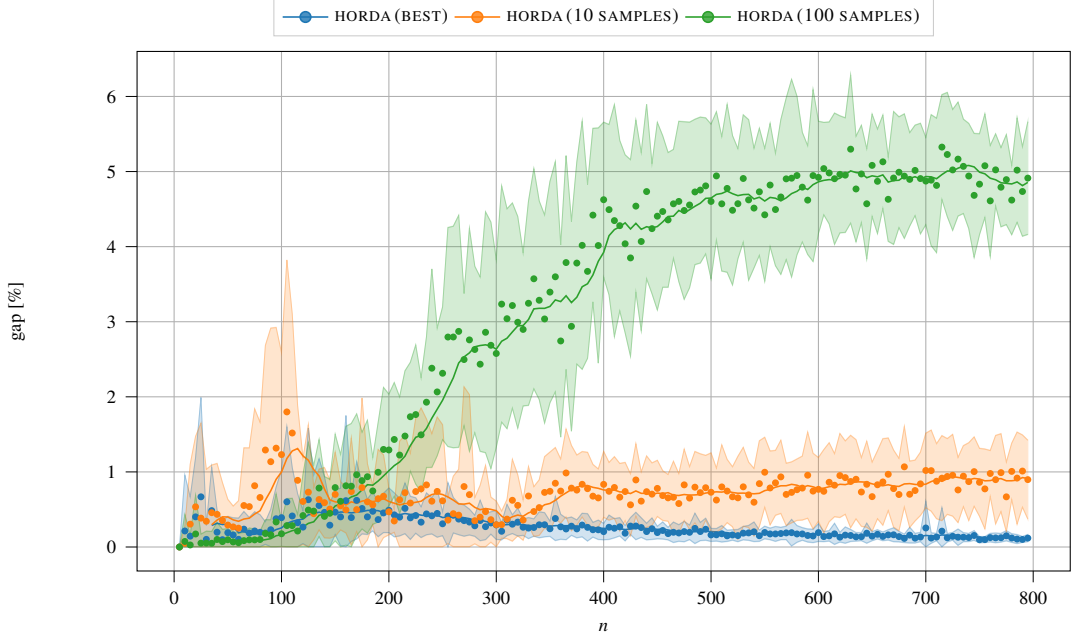}
        \caption{Optimality gap of \DHS{} with a different number of training instances for different $\numjobs{}$.}
        \label{fig:exp:nn:data-size:count}
    \end{figure}

    In this section, we analyze the impact of the neural network on the performance of \DHS{} and on the quality of its solutions. 
    Specifically, we study the impact of the \LSTM capacity, \GRU as an alternative to \LSTM, the number of instances in the training data set, and the size of instances (i.e., $\numjobs$) used for training.
    The experiments listed below assume \DHS{best} as the baseline scenario. This means that every experiment varies one hyperparameter of the neural network while the others are set as in \DHS{best}.
    
    At first, we present the results of \LSTM with different capacities.
    The size of the \LSTM layer impacts the run time of \DHS{} and the ability of the neural network to fit on the training data and their generalization.
    The optimality gaps of \DHS{} with the neural networks having \LSTM capacity $32$, $64$, $128$, $256$ are shown in \cref{fig:exp:nn:capacity-gap}.
    For those capacities, the number of trainable parameters inside the neural network is equal to $4500$, $1.7\powten{4}$, $6.7\powten{4}$, $2.6\powten{5}$, respectively.
    Let us recall that the best model has a capacity equal to $256$, and results show that it performs the best.
    The average ability to fit the training data is decreasing with the decreasing capacity of the \LSTM layer, and thus the optimality gap of \DHS{} grows.
    On the other hand, larger capacities than $256$ would require enormous training data set due to a huge number of training parameters.

    An alternative to the \LSTM layer used in the neural network is the \GRU layer. 
    \GRU is a more restricted architecture (than \LSTM), which can lead to a weaker ability of generalization, whereas the \LSTM architecture can profit from a large number of training samples.
    The comparison of the best model and neural network with the \GRU layer is shown in \cref{fig:exp:nn:capacity-gru}.
    \DHS{GRU}, i.e., the algorithm with the neural network containing the \GRU layer, provides solutions with similar quality in terms of the optimality gap for instances with up to $450$ jobs.
    For the instances with more than $450$ jobs, the optimality gap grows up to $1\%$.
    Therefore, for our purposes is better to use \LSTM, since it is computationally inexpensive to generate large training data sets with the \gensubinstance method.

    The following experiment, illustrated in \cref{fig:exp:nn:data-size:move-n}, evaluates the impact of the number of jobs $\numjobs$ of instances used to generate the training data set.
    The experiment assumes \gensubinstance, where the training used input instances with $50$ - $75$, $75$ - $100$, and $100$ - $125$ jobs.
    For each size of input instances, we generate $20$ instances and generate the training samples that include subproblems as it is described in \cref{sec:sol:data-gen}.
    The best model is trained with input instances of size $75$ - $100$, i.e., \DHS{best}.
    \DHS{} with a neural network trained on the data set generated from instances with $100$ - $125$ jobs has a slowly growing optimality gap.
    For the training on the data set with smaller input instances $50$ - $75$, we observe the peak on the optimality gap at $150$. 
    A similar peak on the gap curve was observed even in other experiments; however, it is interesting how the change of the training instances' size affects the position of the peak and its width.
    Nevertheless, \DHS{best} is superior to either \DHS{$50$-$75$ \n}, and \DHS{$100$-$125$ \n}.
   
    The training data set can be created by \gensubinstance using \edd decomposition, \spt decomposition, or their combination denoted \shorter where the algorithm always selects the decomposition having the smallest $\overline{K^{\circ}}(\jobsset)$.
    Those three possibilities are denoted as \DHS{\edd data set}, \DHS{\spt data set}, and \DHS{best}, respectively.
    The optimality gap of \DHS{} using neural networks trained on data set generated with different decompositions is shown in \cref{fig:exp:nn:data-gen}.
    \DHS{best} has the smallest optimality gap, as it utilizes the same decomposition for the generation of the training data set and the evaluation.
    This underlines the importance of training a neural network on the data with the same distribution as in the evaluation phase achieved by \gensubinstance generator.
    
    A common practice of improving the performance of the neural network is enlarging the training data set; therefore, the following experiment is focused on the effect of its size.
    The neural network was trained on sub-instances generated by \gensubinstance where for each $\numjobs \in \left[75, 100\right]$ we generated $10$, $20$, $100$ input instances.
    In other words, the data sets consist $7.5 \cdot 10^5$, $1.6 \cdot 10^6$, and $7.7 \cdot 10^6$ training samples, respectively.
    The corresponding experiments shown in \cref{fig:exp:nn:data-size:count} are denoted as \DHS{10 samples}, \DHS{best}, and \DHS{100 samples}.
    \DHS{10 samples} has the optimality gap under $1\%$.
    By increasing the training data set's size, we get \DHS{best} and the worst gap decreases under $0.5\%$.
    \DHS{100 samples} has the smallest optimality gap for the instances with up to 125 jobs.
    However, for larger instances, the optimality gap grows rapidly.
    It is important to observe that \DHS{100 samples} has the optimality gap smaller than \DHS{best} on the range of instances from the training data set, i.e., instances with $75-100$ jobs.
    Thus, even though it seems counter-intuitive at first, increasing the training data set the size above a certain critical level can lead to overfitting, resulting in worse prediction performance of the model on the instances that do not lay in a range of the training instances.

\subsection{Parameters of the Scheduling Algorithm}
    \label{sec:exp:horda-params}
    \begin{figure}[ht!]
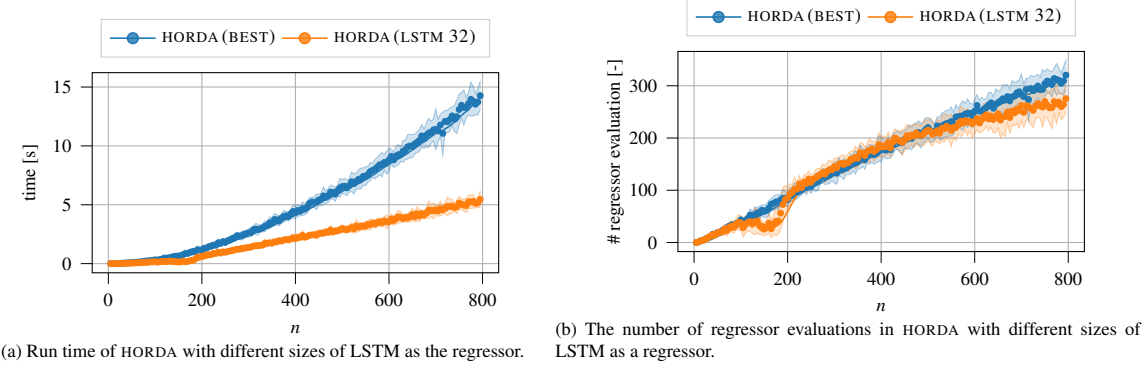

    \centering
        \begin{subfigure}{.5\textwidth}
            \centering
            \captionsetup{skip=0pt}
            \input{img/time_dr_asc_horda_runtime_n800_gbn_sethardest_colinstances_validationavgstdDevPop}
            \caption{Run time of \DHS{} with different sizes of \LSTM as the regressor.}
            \label{fig:exp:dhs:runtime:capacity}
        \end{subfigure}%
        \begin{subfigure}{.5\textwidth}
            \centering
            \captionsetup{skip=0pt}
            \input{img/bonus_regressor_count_sum_dr_asc_horda_runtime}
            \caption{The number of regressor evaluations in \DHS{} with different sizes of \LSTM as a regressor.}
            \label{fig:exp:dhs:shift:eval}
        \end{subfigure}
        \caption{Impact of the \LSTM capacity on the run time and number of regressor calls in \DHS{}.}
    \end{figure}


    \begin{figure}[ht!]
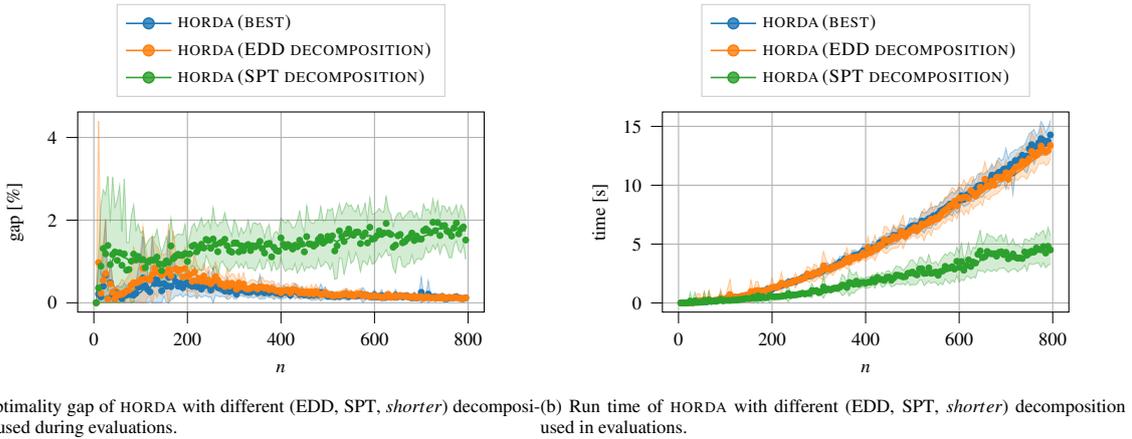

    \centering
        \begin{subfigure}{.5\textwidth}
            \centering
            \input{img/gap_dr_diff_subproblem_eval_n800_gbn_sethardest_colinstances_validationavgstdDevPop}
            \caption{{\updatemb{}Optimality gap of \DHS{} with different (\edd, \spt, \shorter) decompositions used during evaluations.}}
            \label{fig:exp:dhs:decomposition:gap}
        \end{subfigure}%
        \begin{subfigure}{.5\textwidth}
            \centering
            \input{img/time_dr_diff_subproblem_eval_n800_gbn_sethardest_colinstances_validationavgstdDevPop}
            \caption{Run time of \DHS{} with different (\edd, \spt, \shorter) decomposition used in evaluations.}
            \label{fig:exp:dhs:decomposition:runtime}
        \end{subfigure}
        \caption{Impact of the decomposition used in \DHS{} to the optimality gap and run time.}
    \end{figure}
    
    In this section, we study the impact of the \DHS{} parameters on its performance. 
    At first, we analyze how the time of the inference from the neural network affects the run time of \DHS{}.
    Then, we compare results for different decomposition rules used in \DHS{}.

    In the first experiment, we present the impact of the neural network inference time on the run time of \DHS{}.
    For the illustration, we use two neural networks with different inference times.
    The inference time of our neural network is primarily affected by the capacity of the \LSTM layer.
    The faster neural network has a capacity equal to $32$. The slower one has a capacity $256$, which corresponds to \DHS{best} scenario.
    The other parameters are set according to \DHS{best} scenario.
    The run times of \DHS{} assuming capacity $32$ and $256$ (denoted \DHS{LSTM 32} and \DHS{best}, respectively) are shown in \cref{fig:exp:dhs:runtime:capacity}.
    \DHS{LSTM 32} is about three times faster on the instance with 800 jobs than \DHS{best}, and this difference increases with the size of the instance.
    Apart from that, \DHS{LSTM 32} run time slightly drops for instances around 150 jobs. 
    This is caused by a lower number of neural network evaluations, illustrated in \cref{fig:exp:dhs:shift:eval}, where the drop near 150 jobs is obvious even more.
    The lower number of evaluations is connected with a phenomenon we observed in \cref{fig:exp:nn:capacity-gap} that compares the optimality gap for different capacities of \LSTM. 
    In that figure, one can see a sudden deterioration of the solution quality for instances around $\numjobs = 150$.
    The most probable explanation of the correlation between the quality of results and the number of regressor evaluations is that \DHS{} with a poor regressor makes some decisions that lead to subproblems where filtering rules dramatically reduce the candidate set $\positionsset{}$. 
    Due to this, we see fewer regressor evaluations in \cref{fig:exp:dhs:shift:eval}.
    However, these choices are suboptimal and lead to particularly degenerative solutions with poorer quality.
      
    The last experiment presents the impact of the decomposition used in \DHS{} on the quality of the solutions and run times.
    Specifically, it compares the cases when position sets \positionsset{} are generated either by \edd, \spt, or \shorter decomposition.
    \cref{fig:exp:dhs:decomposition:gap} and \cref{fig:exp:dhs:decomposition:runtime} show the solution quality and run time for different decompositions used in \DHS{}.
    \DHS{best} utilizes \shorter decomposition, \DHS{\spt decomposition} utilizes the \spt decomposition, and \DHS{\edd decomposition} utilizes the \edd decomposition.
    The quality of \DHS{\spt decomposition} solutions is inferior to other methods; the results of \DHS{best} and \DHS{\edd decomposition} are similar.
    This is caused by the fact that during the evaluation of \DHS{} with the \shorter decomposition, the \edd decomposition is selected more frequently than the \spt decomposition.
    Since the neural network is trained on the data set generated by \shorter decomposition, the data set contains fewer samples related to \spt decomposition.
    This property can lead to relatively small differences in results between the \edd and \shorter decomposition and is significantly different from the \spt decomposition.
    \DHS{best} provides a better solution than \DHS{\edd decomposition}, mainly for the instances with less than $300$ jobs; for the bigger instances, the difference is negligible.

\section{Conclusion}
\label{sec:conclusion}
To the best of our knowledge, this is one of the first scheduling algorithms where deep learning is successfully used to guide solution-space exploration.
Our approach lies in the synergy between the state-of-the-art operations research method and our neural network. This is opposite to the classical approach in \ML, e.g., of Vinyals et al.~\cite{vinyals2015} with an end-to-end approach for Traveling Salesman Problem.
For the single machine scheduling problem minimizing total tardiness, we show how a neural network can extend standard decomposition techniques. 
Besides, we provide an efficient way to generate the training data set, which is a very time costly operation for combinatorial problems.
The experimental results show that our approach provides near-optimal solutions very quickly and is also able to generalize the acquired knowledge to larger instances without significantly affecting the quality of the solutions.
Our approach has an average gap $0.26\%$ for instances with up to $800$ jobs and outperforms state-of-the-art constructive heuristic \NBR{} with gap $2.14\%$, as well as the decomposition-based heuristic having  gap $1.25\%$.
Moreover, with limited time to 15\jed{s} the state-of-the-art exact algorithms~\cite{Garraffa2018} have an average gap $1.81\%$ and is also dominated by our approach in this scenario.

We believe that the proposed methodology opens new possibilities for the design of efficient heuristics algorithms where the manual tuning of the heuristic is substituted by automatic \ML. Therefore, future research should address other simple scheduling problems that cannot be efficiently decomposed, like \SMTTP. One possibility may be problem $1||\sum w_j T_j$ which is still simple and suitable for neural networks. Another research direction is the generation of the training data set and its efficiency for \nphard{} problems. This paper has shown that there are better ways to generate the training instances; nevertheless, it is tailored to problem \SMTTP.

\section{Acknowledgements}
The authors want to thank Vincent T'Kindt from Université de Tours for providing the source code of \TTBR{} algorithm.

This work was supported by the European Regional Development Fund under the project AI\&Reasoning (reg. no. CZ.02.1.01/0.0/0.0/15\_003/0000466), the Grant Agency of the Czech Republic under the Project GACR 22-31670S, and the EU and the Ministry of Industry and Trade of the Czech Republic under the Project OP PIK CZ.01.1.02/0.0/0.0/20\_321/0024399.

\section*{Appendix}
\label{sec:appendix}

        \begin{table}[thb!]
            \centering
            \caption{List of used notations.
            }
            \begin{tabular}{|c|c|}
                \toprule
            \jobsset & set of jobs (problem instance)\\
            \n & number of jobs\\
            $\proctime{\job}$ & processing time of job $\job$\\
            $\duedate{\job}$ & due date of job $\job$\\
            \sequence{} & permutation of jobs \jobsset\\
            $\objopt{\jobsset}$ & optimal total tardiness of instance \jobsset\\
            $\mathcal{\tardiness}_{\sequence[\position]}(\jobsset)$ & tardiness of job $\sequence[\position]$ in permutation \sequence of \jobsset \\
            $\obj{\jobsset,\sequence}$ & total tardiness of \jobsset under permutation \sequence \\
            $\decabstraction$ & problem decomposition (\edd or \spt)\\
            \jsplit & splitting job in decomposition $\decabstraction$ \\
            \position & possible positions of \jsplit in the schedule\\
            \positionsset{} & set of positions \position of job \jsplit defined by decomposition \decabstraction \\
            $\overline{K^{\circ}}(\jobsset)$ & filtered \positionsset\\
            \precprob{\jobssetrec,\position} & preceding subset of jobs for decomposition \decabstraction, set \jobssetrec and position \position\\
            \follprob{\jobssetrec,\position} & following subset of jobs  for decomposition \decabstraction, set \jobssetrec and position \position\\
            $\jobsset^\prime$ & a subproblem (either $\precprob{\jobssetrec,\position} \subset \jobsset$ or $\follprob{\jobssetrec,\position} \subset \jobsset$)\\
            $Q(\jobsset,\position)$ & optimal total tardiness of \jobsset with splitting job \jsplit at position \position\\
            $\hat{T}(\jobsset)$ & estimated total tardiness of \jobsset\\
            $\hat{Q}(\jobsset,\position)$ & estimated total tardiness of \jobsset with splitting job \jsplit at position \position\\
            \longestjobposition & position of splitting job \jsplit with minimal $\hat{Q}(\jobsset,\position)$\\
            \inp{} & input vector of the neural network\\
            \out{} & output of the neural network\\
            \eddsequence{} & permutation of \jobsset in \edd order\\
            $gap_{EDD}$ & gap of the EDD schedule w.r.t. to the optimal solution\\
            $\maxproc$ & maximal processing time of a job\\
                \bottomrule
            \end{tabular}
            \label{tab:app:notation}
        \end{table}
    
        \begingroup
        \begin{table}[thb!]
        
            \centering
            \caption{List of used abbreviations.}
            \begin{tabular}{|c|c|}
                \toprule
                    \acs{lstm} & Long Short-Term Memory\\
                    \acs{gru} & Gated Recurrent Unit\\
                    \acs{ml} & Machine Learning\\
                    \acs{ga} & Genetic Algorithm\\
                    \acs{dhs} & Heuristic Optimizer using Regression-based Decomposition Algorithm\\
                    \acrotable{ttbr}
                    \acrotable{tsp}
                    \acs{edd} & Earliest Due Date\\
                    \acs{spt} & Shortest Processing Time\\
                    \acs{or}  & Operations Research\\
                    \acrotable{csp}
                    \acs{mdd} & Modified Due Date Rule\\
                \bottomrule
            \end{tabular}
            \label{tab:app:abbr}
        \end{table}
    \endgroup
        
    
\clearpage
\bibliographystyle{acm}
\bibliography{main.bib}




\end{document}